\documentclass[]{siamart250211}

\usepackage{bookmark}
\usepackage{amssymb}
\usepackage{stmaryrd}
\SetSymbolFont{stmry}{bold}{U}{stmry}{m}{n}
\usepackage{array}
\usepackage{multirow}
\usepackage{arydshln}
\usepackage{diagbox}
\usepackage{makecell}
\usepackage{subcaption}
\usepackage{tikz}

\usetikzlibrary{decorations.pathreplacing}
\allowdisplaybreaks

\def \i {{\mathrm{i}}}
\def \ri {{\rm i}}
\def \kr {k_\rho}
\newcommand{\J}[1]{\bs{J}_{#1}}
\newcommand{\bs}[1]{\boldsymbol{#1}} 

\newcommand{\bd}[1]{\boldsymbol{#1}}

\usepackage{newunicodechar}
\newunicodechar{−}{\ensuremath{-}}

\usepackage{algorithm}
\usepackage{algorithmicx}
\usepackage[noend]{algpseudocode}
\newsiamremark{remark}{Remark}
\newsiamremark{example}{Example}
\usepackage{booktabs}

\title{Fast Multipole method for Maxwell's equations in layered media}
\author{ Heng Yuan\thanks{LCSM(MOE), School of Mathematics and Statistics, Hunan Normal University, Changsha, Hunan, 410081, P. R. China.}
        \and Bo Wang\thanks{Joint first author. LCSM(MOE), School of Mathematics and Statistics, Hunan Normal University, Changsha, Hunan, 410081, P. R. China. This author acknowledges financial support provided by NSF of China (grant 11771137, 12022104), and by the Construct Program of the Key Discipline in Hunan Province.}
        \and Wenzhong Zhang\thanks{Corresponding author, Suzhou Institute for Advanced Research, University of Science and Technology of China, Suzhou, Jiangsu, 215127, P. R. China. ({\tt wenzhong@ustc.edu.cn}). This author acknowledges financial support provided by NSF of China (grant 12201603). \textit{Submitted to SIAM Journal on Scientific Computing, July 25, 2025.} }
        \and Wei Cai\thanks{Department of Mathematics, Southern Methodist University, Dallas, TX 75275, USA.}
        }%

\begin{document}
\maketitle

\begin{abstract}
We present a fast multipole method (FMM) for solving Maxwell's equations in three-dimensional (3-D) layered media, based on the magnetic vector potential $\bs A$ under the Lorenz gauge, to derive the layered dyadic Green’s function. The dyadic Green's function is represented using three scalar Helmholtz layered Green’s functions, with all interface-induced reaction field components expressed through a unified integral representation. By introducing equivalent polarization images for sources and effective locations for targets to reflect the actual transmission distance of different reaction field components, multiple expansions (MEs) and local expansions (LEs) are derived for the far-field governed by actual transmission distance. To further enhance computational efficiency and numerical stability, we employ a Chebyshev polynomial expansion of the associated Legendre functions to speed up the calculation of multipole-to-local (M2L) expansion translations. Finally, leveraging the FMM framework of the Helmholtz equation in 3-D layered media \cite{bo2019hfmm},  we develop a FMM for the dyadic Green’s function of Maxwell's equations in layered media. Numerical experiments demonstrate the $\mathcal O(N\log N)$-complexity of the resulting FMM method, and rapid convergence for interactions of low-frequency electromagnetic wave sources in 3-D layered media.
\end{abstract}



\section{Introduction}

The fast multipole method (FMM) \cite{GreenGardRokhlin1987, GreenGardRokhlin1997}, originally introduced to accelerate the evaluation of pairwise interactions in gravitational and electrostatic $N$-body systems, has become a revolutionary numerical technique in computational electromagnetics (CEM) over the past three decades. The key contribution is its ability to reduce the computational complexity of evaluating long-range interactions—from the direct $\mathcal O(N^2)$ scaling to 
$\mathcal O(N)$ or $\mathcal O(N\log N)$—by exploiting hierarchical domain decomposition and analytic expansions of the underlying kernel functions. 
FMM has changed the field of computational electromagnetics by enabling the solution of large-scale problems in $O(N)$-linear scaling of computing time. In free space, under the Lorenz gauge, the dyadic Green's functions for Maxwell's equations can be constructed by applying a differential operator to the scalar Green's function of the Helmholtz equation. This enables the derivation of FMM for Maxwell's equations from the Helmholtz FMM framework. J.-M Song and W. C. Chow proposed  the MLFMA \cite{song1995Multilevel, chew2001fast}, an extension of the FMM for vector wave equations, to efficiently handle large-scale electromagnetic scattering problems \cite{chew2001fast,tong2020nystrom}. This approach is feasible to simulate structures with millions of unknowns.

Over the past two decades, sustained efforts have been devoted to develop fast algorithms to solve electromagnetic scattering problems in layered media (cf. \cite{hu2000fast,pan2004fast,chew2005fast,chen2018accurate,chen2017pml,lu2019numerical,cho2015robust,gueuning2022inhomogeneous,gao2022wave,arrieta2022windowed}), driven by their critical applications in very large-scale integrated (VLSI) circuit simulation \cite{kao2001parasitic,weeks2003calculation,chew2005toward,karsilayan2010full,sharma2018complete,sharma2020slim}, geophysics (cf. \cite{he2000method,okhmatovski2024electromagnetic}), medical imaging (cf. \cite{xu2009image,weber2006modulated,chen2017electromagnetic}), and etc. However, the presence of material interfaces poses significant challenges because of the reflection and transmission of electromagnetic waves across layers. A straightforward method employs the free-space dyadic Green's functions and introduces additional unknowns on material interfaces in order to enforce transmission conditions on interfaces. Such an approach significantly increases the size of the resulting linear system, particularly when the number of layers is large. Utilizing layered dyadic Green's functions \cite{TETM2002,chew2006matrix} satisfying the transmission conditions on the material interfaces is a more natural approach, avoiding the involvement of additional unknowns on those interfaces. However, this creates a challenge for an efficient and accurate computation of the source interactions governed by the layered dyadic Green's function, necessitating the development of specialized fast algorithms.

In a series of recent works, we have systematically developed a unified framework for constructing FMMs for the 3-D Helmholtz equation \cite{bo2019hfmm} and other important scalar physical equations in layered media \cite{zhang2020, bo2021laplace, bo2021PoissonBoltzmann,zhang2022exponential}. In this framework, the interaction among sources in layered media is decomposed into the free-space components identical to the free-space problem, as well as the reaction field components incited by the layered structure. The layered FMM approach integrates conventional free-space FMM with newly developed fast algorithms for the reaction field components, achieving computational efficiency comparable to that of classic free-space versions.

In this work, we extend the works of layered FMM of scalar equations to the $3 \times 3$ dyadic Green's functions associated with Maxwell's equations in layered media, aiming at enabling efficient integral equation solvers for electromagnetic scattering problems in stratified environments. Building upon the integral formulations introduced in \cite{bo2022maxwellDGF}, we derive a unified representation of layered dyadic Green's functions that reveals $5$ distinct categories of angular dependence in the Fourier spectral domain, so that the far-field expansions are neatly introduced without the need of derivatives on LE basis functions.
We also improve the method of equivalent polarization coordinates from our previous implementation \cite{bo2019hfmm,zhang2020} by more carefully evaluating the vertical transmission distance of reaction field components, and employing polarization coordinates for sources and effective locations for targets.
Based on the concept of polarization coordinates and effective locations, we construct the corresponding far-field expansions and develop efficient shifting and translation operators for the five categories of Sommerfeld-type integrals using the extended Funk–Hecke identity (cf. \cite{bo2019hfmm}).
To further accelerate the precomputation of the multipole-to-local (M2L) translation matrices, we incorporate a Chebyshev polynomial expansion for the products of two associated Legendre functions.
With these novel far-field approximation formulas, a complete FMM is established for evaluating the vector potential induced by directed Hertz dipole current sources in layered media.
Numerical experiments in two-layers and three-layers configurations confirm the $\mathcal{O}(N\log N)$ complexity and spectral convergence with respect to the truncation parameter $p$ of the proposed method. Similarly to the FMM developed for the Helmholtz equation in 3-D layered media, the reaction field components incur significantly lower computational cost than their free-space counterparts, due to the general separation of equivalent polarization coordinates of sources and effective locations of targets. As a result, the overall computational complexity of the proposed method remains comparable to that of free-space FMMs, as long as the number of layers in the media is not large.

The rest of the paper is organized as follows.
\Cref{sect-freespace} provides a brief review of the FMM for Maxwell's equations in free space.
\Cref{sect-layeredmedia} elaborates the formulation and implementation of the layered media FMM for Maxwell's equations, including the derivation of the far-field approximation for the reaction field components, as well as various techniques to improve the efficiency in the evaluation of M2L translations. In \Cref{sect-num}, numerical examples are provided to validate the accuracy and efficiency of the proposed method. The final conclusion is given in \Cref{sect-conclusion}.

\section{A review of the FMM for Maxwell's equations in free space}\label{sect-freespace}
In this section we briefly review the multipole and local expansions and their shifting and translation operators of the Maxwell's equations in the free space.

We assume a time dependence $e^{\ri\omega t}$ in Maxwell's equations throughout this paper, where $\omega$ is the angular frequency in time.
The interaction between a target particle $\bs{r} = (x, y, z) \in \mathbb{R}^3$ and a source particle $\bs{r}' = (x', y', z') \in \mathbb{R}^3$ is discussed for simplified illustration.
The electric field and magnetic field dyadic Green's functions $\bs{G}_E^f(\bs{r};\bs{r}')$ and $\bs{G}_H^f(\bs{r};\bs{r}')$ of the Maxwell's equations in the free space are defined using a $3 \times 3$ potential tensor $\bs{G}_A^f(\bs{r};\bs{r}')$ by 
\begin{equation}\label{eq-GEGH-GA}
\bs{G}_{E}^f = -\i \omega \left(\bs{I} + \frac{\nabla\nabla}{k^2}\right)\bs{G}_{A}^f,\quad \bs{G}_{H}^f = \frac{1}{\mu} \nabla \times \bs{G}_{A}^f,
\end{equation}
where $\epsilon_0$, $\mu_0$ are the dielectric constant and magnetic permeability in vacuum, respectively, and $ k=\omega\sqrt{\epsilon_0\mu_0}$ is the wave number in vacuum. The potential tensor 
\begin{equation*}
	\bs{G}_{\bs A}^f(\bs{r};\bs{r}') = -\frac{1}{\i \omega} g^f(\bs{x};\bs{x}') \bs{I},\quad g^f(\bs{r};\bs{r}') = \frac{e^{\ri k|\bs{r}-\bs{r}'|}}{4\pi|\bs{r}-\bs{r}'|}=\dfrac{\ri  k}{4\pi}h_0^{(1)}( k|\bs{r}-\bs{r}'|)
\end{equation*} 
is the solution to the vector Helmholtz equation
\begin{equation}\label{eq-GA-eqn}
\nabla^2 \bs{G}_A^f + k^2 \bs{G}_A^f = \frac{1}{\i \omega} \delta(\bs{r}-\bs{r}')\bs{I}.
\end{equation}
In short, the dyadic Green's functions are given by
\begin{align}\label{GEandGH}
	\bs{G}_{\bs E}^f(\bs{r};\bs{r}') = \left(\bs{I} + \frac{\nabla\nabla}{k^2} \right)g^f(\bs r, \bs r') ,
	\quad \bs{G}_{\bs H}^f(\bs{r};\bs{r}') = \frac{\ri}{\omega\mu_0} \nabla \times (g^f(\bs r, \bs r') \bs I),
\end{align}
respectively.

For the scalar Green's function $g^f(\bs r;\bs r')$ of Helmholtz equation with wave number $k$, we have a multipole expansion with respect to a (source) center $\bs{r}_c^s$
\begin{equation}\label{freespace3dmulexpscaled}
	g^f(\bs r; \bs r')=\ri k\sum\limits_{n=0}^{\infty}\sum\limits_{m=-n}^n M_{nm}(\bs r_c^s)\frac{h_n^{(1)}(kr_s)}{h_n^{(1)}(kS)}Y_n^m(\theta_s,\varphi_s),
\end{equation}
and a local expansion with respect to a (target) center $\bs{r}_c^t$
\begin{equation}\label{freespace3dlocexpscaled}
	g^f(\bs r; \bs r')=\ri k\sum\limits_{n=0}^{\infty}\sum\limits_{m=-n}^n L_{nm}(\bs r_c^t)h_n^{(1)}(kS')j_n(kr_t)Y_n^m(\theta_t,\varphi_t),
\end{equation}
where the scalar multipole expansion and local expansion coefficients are given by
\begin{equation}\label{melecoeffreescaled}
	 M_{nm}(\bs r_c^t)=h_n^{(1)}(kS)j_n(kr_s')\overline{Y_n^{m}(\theta'_s,\varphi_s')},\quad  L_{nm}(\bs r_c^t)=\frac{h^{(1)}_n(kr_t')}{h^{(1)}_n(kS')}\overline{Y_n^{m}(\theta_t',\varphi_t')},
\end{equation}
respectively, where $\bs{r}^s_c = (x^s_c,y^s_c,z^s_c)$ is the source center close to $\bs{r}'$, $\bs{r}^t_c = (x^t_c,y^t_c,z^t_c)$ is the target center close to $\bs{r}$, $(r_s,\theta_s,\varphi_s)$, $(r_t,\theta_t,\varphi_t)$ are the spherical coordinates of $\bs{r}-\bs{r}^s_c$ and $\bs{r}-\bs{r}^t_c$, respectively, and $(r_s',\theta_s',\varphi_s')$, $(r_t',\theta_t',\varphi_t')$ are the spherical coordinates of $\bs{r}'-\bs{r}^s_c$ and $\bs{r}'-\bs{r}^t_c$, respectively. Here, scaling factors $h_n^{(1)}(kS)$ and $h_n^{(1)}(kS')$ with characteristic lengths $S$ and $S'$ are introduced to avoid possible overflow and underflow in the numerical implementation of the FMM. The characteristic lengths are chosen as sizes of the boxes in the source and target tree from the hierarchical structure of FMM, respectively.
In our tests, we find these scaling factors superior to the power scaling $S^n$ (see e.g. \cite{dashmm16}) in terms of numerical stability.

The far-field expansions of the dyadic Green's functions of Maxwell's equation are straightforwardly derived by applying the tensor differential operators in \cref{GEandGH} to the expansions of the scalar Green's function $g^f(\bs r; \bs r')$. We use $\bs{G}_{\bs E}^f$ as an example.
By merging the expansion \cref{freespace3dmulexpscaled} into \cref{GEandGH}, the multipole expansion (ME) of the electric field dyadic Green's function $\bs{G}_{\bs E}^f$ is given as
\begin{equation}\label{ME_Maxwell_GE}
	\bs{G}_{\bs E}^f(\bs{r};\bs{r}')=\ri k\sum\limits_{n=0}^{\infty}\sum\limits_{m=-n}^n M_{nm}(\bs r_c^s)\left(\bs{I} + \frac{\nabla\nabla}{k^2} \right)\frac{h_n^{(1)}(kr_s)}{h_n^{(1)}(kS)}Y_n^m(\theta_s,\varphi_s),
\end{equation}
and similarly the local expansion (LE) as
\begin{equation}\label{LE_Maxwell_GE}
	\bs{G}_{\bs E}^f(\bs{r};\bs{r}')=\ri k\sum\limits_{n=0}^{\infty}\sum\limits_{m=-n}^n L_{nm}(\bs r_c^t)\left(\bs{I} + \frac{\nabla\nabla}{k^2} \right)h_n^{(1)}(kS')j_n(kr_t)Y_n^m(\theta_t,\varphi_t).
\end{equation}
The translations from multipole expansion to local expansion (M2L), multipole expansion to multipole expansion (M2M) and local expansion to local expansion (L2L) are consistent with those for the scalar Helmholtz equation, which are given as follows:
\begin{align}
	 L_{nm}(\bs r_c^t)&=\sum_{\nu=0}^{\infty}\sum_{\mu= -\nu}^{\nu}\dfrac{S_{n\nu}^{m\mu}(\bs r_c^t-\bs r_c^s)}{h_{\nu}^{(1)}(kS)h_n^{(1)}(kS')} M_{\nu\mu}(\bs r_c^s), \label{maxwell_M2L_freespace} \\
	 M_{nm}(\tilde{\bs r}_c^s)&=\sum_{\nu=0}^{\infty}\sum_{\mu= -\nu}^{\nu}\overline{\widehat{ S}_{n\nu}^{m\mu}(\bs r_c^s-\tilde{\bs r}_c^s)}\dfrac{h_{n}^{(1)}(k{\widetilde{S}})}{h_{\nu}^{(1)}(kS)} M_{\nu\mu}(\bs r_c^s),\label{maxwell_M2M_freespace}\\
	L_{nm}(\tilde{\bs r}_c^t)&=\sum_{\nu=0}^{\infty}\sum_{\mu= -\nu}^{\nu}\widehat{S}_{\nu n}^{\mu m}(\tilde{\bs r}_c^t-\bs r_c^t)\dfrac{h_{\nu}^{(1)}(kS')}{h_n^{(1)}(k\widetilde{S}')} L_{\nu\mu}(\bs r_c^t),\label{maxwell_L2L_freespace}
\end{align}
where $\{M_{nm}(\tilde{\bs r}_c^s)\}$ and $\{L_{nm}(\tilde{\bs r}_c^t)\}$ are the multipole and local expansion coefficients with respect to new source center $\tilde{\bs r}_c^s$ and target center $\tilde{\bs r}_c^t$, respectively, $(S_{n\nu}^{m\mu})$, $(\widehat S_{n\nu}^{m\mu})$ are separation matrices used in the addition theorem of wave functions \cite{martin2006multiple}, and $\widetilde{S},\widetilde{S}'$ are characteristic lengths associated to the new centers. In the implementation of the FMM, the M2M translation proceeds from bottom to top in the FMM tree structure. Specifically, it converts the multipole expansions (MEs) of the eight child boxes at a lower level into the ME of their parent box at the next higher level. Conversely, the L2L  translation goes from top to bottom, transforming the local expansions (LEs) of a parent box into the LEs of its eight child boxes at the lower level.
As a result, in the FMM implementation, we have $\widetilde{S} = 2S$ and $\widetilde{S}' = S'/2$.

\begin{remark}
In the multipole expansion \cref{ME_Maxwell_GE} and local expansion \cref{LE_Maxwell_GE}, the tensor differential operators are evaluated upon the scaled special functions
in the implementation of FMM.
They can be recursively calculated at $\mathcal{O}(p^2)$ cost.
We emphasize that, later in our proposed method handling the layered counterpart, such Hessian matrices are not necessarily evaluated. The far-field approximations are derived directly for ${\bs G}_{\bs E}$ instead of the potential ${\bs G}_{\bs A}$. Thus, the far-field approximations \cref{maxwell_ME_layered,maxwell_LE_layered} derived for the reaction field components have higher order of convergence than \cref{ME_Maxwell_GE,LE_Maxwell_GE}.
\end{remark}

\section{FMM for Maxwell's equations in layered media}\label{sect-layeredmedia}
In this section, we introduce a concise formulation of the dyadic Green's functions of Maxwell's equations in layered media using a matrix basis that separates angular dependence. Then, we propose the far-field expansions including MEs, LEs and their translations based on the formulation.
Finally, we discuss the implementation of the FMM, especially techniques that further improve the performance of M2L evaluation.

\subsection{The dyadic Green's functions of Maxwell's equations in layered media}\label{sect-dyadicGF}
We start from a quick review of \cite{bo2022maxwellDGF}, where in layered media the Maxwell's equations are decomposed into scalar Helmholtz equations.
The approach is consistent with the formulation derived in \cite{TETM2002}.

Consider a horizontally layered medium with $L+1$ layers indexed by $0, \cdots, L$ from top to bottom, respectively. Let the interface between layer $\ell$ and layer $\ell + 1$ be given by the plane $z=d_\ell$, $0 \le \ell \le L-1$, where $d_0 > \cdots > d_{L-1}$, and the medium in each layer is homogeneous with permittivity $\varepsilon_\ell$ and permeability $\mu_\ell$, $\ell=0,\cdots,L$, respectively.
Define the wave numbers in the $\ell$-th layer by
\begin{equation}\label{wavenumber_layer}
	k_\ell = \omega\sqrt{ \varepsilon_\ell \mu_\ell},\quad \ell=0,\cdots,L.
\end{equation} 
Applying the 2-D Fourier transform
\begin{equation}\label{eq-2DFT}
f(x,y) = \frac{1}{4\pi^2} \iint_{\mathbb{R}^2} e^{\i k_x (x - x') + \i k_y(y - y')} \widehat{f}(k_x, k_y) dk_x dk_y,
\end{equation}
a matrix form of the dyadic Green's functions $\widehat{\bd{G}}_{\bs E}$, $\widehat{\bd{G}}_{\bs H}$ of Maxwell's equations with directed Hertz dipole current source in the layered medium can be obtained in the frequency domain.
Namely, in the frequency domain with $(k_x, k_y)$ coordinates,
\begin{equation}\label{eq-GEGH-b-latest}
    \begin{aligned}
        \widehat{\bd{G}}_{\bs E} &= -\frac{\ri \omega}{k_{\ell}^2} \Big( k_{\ell}^2 \phi_{\ell} \J{1} + \mu_{\ell} \kr^2 \psi_{\ell} \J{2} + \mu_{\ell}\partial_z \psi_{\ell} \J{3} + \mu_{\ell} \tilde\psi_{\ell} \J{4} + \frac{k_{\ell}^2\phi_{\ell } + \mu_{\ell}\partial_z \tilde\psi_{\ell}}{\kr^2}\J{5} \Big)+\bs\delta, \\
        \widehat{\bd{G}}_{\bs H} &= \frac{1}{\mu_{\ell}}\Big( \phi_{\ell} \J{6} + \mu_{\ell} \psi_{\ell} \J{7} + \frac{\partial_z \phi_{\ell }- \mu_{\ell}\tilde\psi_{\ell}}{\kr^2}\J{8} - \partial_z \phi_{\ell } \J{9} \Big)
    \end{aligned}
\end{equation}
for $d_{\ell}<z<d_{\ell-1}$, where
$(\kr, \alpha)$ is the polar coordinates of $(k_x, k_y)$, i.e.
\begin{equation}\label{eq-kr-alpha}
	k_x = \kr \cos \alpha, \quad k_y = \kr \sin \alpha, \quad \alpha \in [0, 2\pi),
\end{equation}
the $3 \times 3$ matrices $\J{1},\cdots, \J{9}$ defined as
\begin{equation}\label{mat-J}
	\begin{array}{lll}
		\J{1}  = \begin{bmatrix}
			1 & & \\
			& 1 & \\
			& & 0
		\end{bmatrix}, &
		\J{2}  = \begin{bmatrix}
			0 & & \\
			& 0 & \\
			& & 1
		\end{bmatrix}, &
		\J{3}  = \begin{bmatrix}
			0 & 0 & \i k_x \\
			0 & 0 & \i k_y \\
			0 & 0 & 0
		\end{bmatrix}, \\[15pt]
		\J{4}  = \begin{bmatrix}
			0 & 0 & 0 \\
			0 & 0 & 0 \\
			\i k_x & \i k_y & 0
		\end{bmatrix}, &
		\J{5}  = \begin{bmatrix}
			-k_x^2 & -k_x k_y & 0 \\
			-k_x k_y & -k_y^2 & 0 \\
			0 & 0 & 0
		\end{bmatrix}, &
		\J{6}  = \begin{bmatrix}
			0 & 0 & 0 \\
			0 & 0 & 0 \\
			-\i k_y & \i k_x & 0
		\end{bmatrix}, \\[15pt]
		\J{7}  = \begin{bmatrix}
			0 & 0 & \i k_y \\
			0 & 0 & -\i k_x \\
			0 & 0 & 0
		\end{bmatrix}, &
		\J{8}  = \begin{bmatrix}
			k_x k_y & k_y^2 & 0 \\
			-k_x^2 & -k_x k_y & 0 \\
			0 & 0 & 0
		\end{bmatrix}, &
		\J{9}  = \begin{bmatrix}
			0 & 1 & 0 \\
			-1 & 0 & 0 \\
			0 & 0 & 0
		\end{bmatrix}
	\end{array}
\end{equation}
form a basis, and
\begin{equation*}\label{delta_Ge_hat}
	\bs\delta=\dfrac{\delta_{\ell\ell'}}{k_{\ell}^2}\left[\dfrac{\J{5}}{k_{\rho}^2}-\J{2}\right]\delta(z-z')
\end{equation*}
with $\delta_{\ell\ell'}$ being the Kronecker symbol. The coefficients $\phi_{\ell}, \psi_{\ell}, \tilde\psi_{\ell}$ in \cref{eq-GEGH-b-latest} are functions defined in the $\ell$-th layer. Define piecewise functions
\begin{equation*}
  \phi(k_{\rho}, z, z')=\phi_{\ell}(k_{\rho}, z, z'),\quad \psi(k_{\rho}, z, z')=\psi_{\ell}(k_{\rho}, z, z'), \quad \tilde\psi(k_{\rho}, z, z')=\tilde\psi_{\ell}(k_{\rho}, z, z')
\end{equation*}
for $d_{\ell}<z<d_{\ell-1}$. Then, the functions $\phi$ and $\psi$ are indeed the Green's functions for Helmholtz equations in the layered media, respectively. More precisely, they are the solutions of the following interface problems
\begin{equation}\label{helmholtzb1}
	\begin{cases}
		\displaystyle\partial_{zz} \phi(k_{\rho}, z, z') + k_{\ell z}^2 \phi(k_{\rho}, z, z')=-\frac{\i}{\omega}\delta(z-z'),\quad d_{\ell}<z<d_{\ell-1},\\
		\displaystyle\llbracket \phi\rrbracket=0,\quad \Big\llbracket \frac{1}{\mu}\partial_z\phi\Big\rrbracket=0,
	\end{cases}
\end{equation}
\begin{equation}\label{helmholtzb2}
	\begin{cases}
		\displaystyle	\partial_{zz} \psi(k_{\rho}, z, z') + k_{\ell z}^2 \psi(k_{\rho}, z, z')=-\frac{\i}{\mu\omega}\delta(z-z'),\quad d_{\ell}<z<d_{\ell-1},\\
		\displaystyle	\llbracket \psi\rrbracket=0,\quad \Big\llbracket \frac{1}{\varepsilon}\partial_z\psi\Big\rrbracket=0,
	\end{cases}
\end{equation}
where
\begin{equation*}
k_{\ell z}:=\sqrt{k_{\ell}^2-k_{\rho}^2},\quad \ell=0,1,\cdots,L,
\end{equation*}
with branch cuts $\Im{(k_{\ell z})}\ge0$, $\llbracket \cdot \rrbracket$ denotes the jump across interfaces $\{z=d_{\ell}\}_{\ell=0}^{L-1}$.
The function $\tilde\psi$ is associated with $\psi$ via
\begin{equation*}
    \tilde\psi(\kr, z, z') = - \partial_{z'}\psi(\kr, z, z').
\end{equation*}
An important advantage of using the above formulation \cref{eq-GEGH-b-latest} is that the coefficients of matrices $\J{1}, \cdots, \J{9}$ are guaranteed to be rotationally symmetric functions regardless of the number of layers, i.e. they are functions of $\kr$ without dependence on the polar angle $\alpha$. This property enables us to derive uniform far-field approximation formulas for the reaction field components in the next subsection.

Solutions to the scalar Green's functions of Helmholtz equations in layered media \cref{helmholtzb1,helmholtzb2} can be found analytically in the frequency domain using classic iterative methods with $\mathcal{O}(L)$ computational cost on each $\kr$, see e.g. \cite{Chew1999inhomogenous,bo2020tefmm}.
In this paper, we adopt the concept of reaction field decomposition from our previous works. Namely, the dyadic Green's functions are interpreted as the free-space dyadic Green's function in the $\ell'$-th layer of source particle, as well as the \emph{reaction field} incited by the presence of other layers. The reaction field is further decomposed into $4$ components regarding different vertical field propagation directions by the target and source.
Specifically,
\begin{equation}\label{threedensity}
\begin{split}
    \phi_{\ell}(k_{\rho}, z, z')&= \delta_{\ell\ell'} \phi^f(k_{\rho}, z, z')-\frac{1}{2\omega k_{\ell' z}}\sum_{\ast,\star=\uparrow,\downarrow}\phi_{\ell\ell'}^{\ast\star}(\kr){Z}_{\ell\ell'}^{\ast\star}(k_{\rho}, z,z'),\\
    \psi_{\ell}(k_{\rho}, z, z')&= \delta_{\ell\ell'} \psi^f(k_{\rho}, z, z')-\frac{1}{2\omega\mu_{\ell'} k_{\ell' z}}\sum_{\ast,\star=\uparrow,\downarrow}\psi_{\ell\ell'}^{\ast\star}(\kr){Z}_{\ell\ell'}^{\ast\star}(k_{\rho}, z,z'),\\
    \tilde\psi_{\ell}(k_{\rho}, z, z')&=-\delta_{\ell\ell'} \partial_{z'}\psi^f(k_{\rho}, z, z')+\frac{\ri }{2\omega \mu_{\ell'}}\sum_{\ast,\star=\uparrow,\downarrow}s_{\ell'\ell'}^{\star\star}\psi_{\ell\ell'}^{\ast\star}(\kr){Z}_{\ell\ell'}^{\ast\star}(k_{\rho}, z,z'),
\end{split}
\end{equation}
for $\bs r$ in the $\ell$-th layer and $\bs r'$ in the $\ell'$-th layer.
We omit the formula of functions $\phi^f$ and $\psi^f$, as they are free-space components and won't be touched in the implementation.
The $Z^{\ast\star}$-functions 
\begin{equation}\label{zexponential}
	\begin{split}
    {Z}_{\ell\ell'}^{\uparrow\uparrow}(k_{\rho}, z,z')=&\begin{cases}
			\displaystyle e^{\ri (k_{\ell z}z-k_{\ell' z}z')},\quad \ell<\ell'\\
			\displaystyle e^{\ri (k_{\ell' z}\tau_{\ell'-1}(z')-k_{\ell z}\tau_{\ell}(z))},\quad\ell\geq \ell'
		\end{cases}\\
		{Z}_{\ell\ell'}^{\uparrow\downarrow}(k_{\rho}, z,z')=&\begin{cases}
			\displaystyle e^{\ri (k_{\ell z}z-k_{\ell' z}\tau_{\ell'}(z'))},\quad \ell\leq\ell'\\
			\displaystyle e^{\ri (k_{\ell' z}z'-k_{\ell z}\tau_{\ell }(z))},\quad \ell>\ell'
		\end{cases}\\
        {Z}_{\ell\ell'}^{\downarrow\uparrow}(k_{\rho}, z,z')=&\begin{cases}
			\displaystyle e^{\ri (k_{\ell z}\tau_{\ell-1}(z)-k_{\ell' z}z')},\quad\ell<\ell',\\
			\displaystyle e^{\ri(k_{\ell' z}\tau_{\ell'-1}(z')-k_{\ell z}z)},\quad\ell\geq \ell',
		\end{cases}\\
		{Z}_{\ell\ell'}^{\downarrow\downarrow}(k_{\rho}, z,z')=&\begin{cases}
			\displaystyle e^{\ri (k_{\ell z}\tau_{\ell-1}(z)-k_{\ell' z}\tau_{\ell'}(z'))},\quad\ell\leq \ell'\\
			e^{\ri (k_{\ell' z}z'-k_{\ell z}z)},\quad \ell>\ell',
		\end{cases}
	\end{split}
\end{equation}
are exponential functions indicating the different cases of vertical field propagation directions, where
\begin{equation}
\label{z-reflection}
    \tau_{\ell}(z)=2d_{\ell}-z
\end{equation}
and
\begin{equation}
\label{eq_signs}
    s_{\ell\ell'}^{\ast\downarrow}=\begin{cases}
        \displaystyle 1,& 	\ell\leq \ell',\\
        \displaystyle -1,& \ell>\ell',\\
    \end{cases},\quad s_{\ell\ell'}^{\ast\uparrow}=\begin{cases}
        \displaystyle 1,& 	\ell< \ell',\\
        \displaystyle -1,& \ell\geq \ell',\\
    \end{cases}\quad \ast=\uparrow,\downarrow
\end{equation}
are signs depending on $\ell,\ell'$.
The density functions $\phi^{\ast\star}_{\ell\ell'}(k_{\rho})$, $\psi_{\ell\ell'}^{\ast\star}(k_{\rho})$ depend only on the layered structure and can be calculated using recurrence formulas in a standard manner \cite{bo2020tefmm, Chew1999inhomogenous}.
The involvement of the $\tau(\cdot)$ notation, corresponding to the reflection with respect to nearby interface planes, ensures that the $Z^{\ast\star}$ functions are guaranteed to decay exponentially as $\kr \to +\infty$. The exponential decay rates are indeed determined by the vertical \emph{transmission distance} in the propagation of each reaction field component, which will be further discussed in the next subsection.
\begin{remark}
It should be noted that the formulations in \eqref{zexponential} are not unique as we could move exponential terms $e^{\ri k_{\ell z} d}$ or $e^{\ri k_{\ell'z}d}$ in/out the density functions $\phi_{\ell\ell}^{\ast\star}(k_{\rho})$, $\psi_{\ell\ell}^{\ast\star}(k_{\rho})$. However, formulations in \eqref{zexponential} are selected due to the minimal vertical \emph{transmission distance} in the propagation of each reaction field component.
\end{remark}

Let us focus on the electric field Green's function. The formulation \cref{eq-GEGH-b-latest} and \cref{threedensity} imply
\begin{equation}\label{componets}
	\widehat{\bd{G}}_{\bs E}=\delta_{\ell\ell'}\widehat{\bd{G}}_{\bs E}^{f}+\widehat{\bd{G}}_{\ell\ell'}^{\uparrow\uparrow}+\widehat{\bd{G}}_{\ell\ell'}^{\uparrow\downarrow}+\widehat{\bd{G}}_{\ell\ell'}^{\downarrow\uparrow}+\widehat{\bd{G}}_{\ell\ell'}^{\downarrow\downarrow}
\end{equation}
while $\widehat{\bd{G}}_{\bs E}^{f}$ is the Fourier transform of $\bd{G}_{\bs E}^{f}$ given in \cref{GEandGH}, and the reaction field components
\begin{equation}\label{reactioncomponentspectral}
    \widehat{\bd{G}}_{\ell\ell'}^{\ast\star}(k_x, k_y, z, z')=\i\frac{{Z}_{\ell\ell'}^{\ast\star}(k_{\rho}, z,z')}{2k_{\ell'z}}\bs \Theta_{\ell\ell'}^{\ast\star}(k_x, k_y),\quad \ast,\star \in \{\uparrow, \downarrow\},
\end{equation}
where, depending on the field propagation directions,
\begin{equation}\label{Theta}
\begin{split}
    \bs \Theta_{\ell\ell'}^{\uparrow\downarrow} &= \phi_{\ell\ell' }^{\uparrow\downarrow} \Big(\J{1}+ \frac{1}{\kr^2}\J{5}\Big)+ \frac{\mu_{\ell}\psi_{\ell\ell'}^{\uparrow\downarrow} }{\mu_{\ell'}k_{\ell}^2} \Big(\kr^2\J{2}+\ri k_{\ell z}\J{3}-  \ri k_{\ell' z}\J{4} + \frac{k_{\ell z} k_{\ell' z}}{\kr^2}  \J{5}\Big), \\
    \bs \Theta_{\ell\ell'}^{\uparrow\uparrow} &= \phi_{\ell\ell'}^{\uparrow\uparrow}\Big(\J{1}+ \frac{1}{\kr^2}\J{5}\Big) + \frac{\mu_{\ell} \psi_{\ell\ell'}^{\uparrow\uparrow}}{\mu_{\ell'}k_{\ell}^2}  \Big(\kr^2\J{2}+\ri k_{\ell z}\J{3}+\ri k_{\ell' z}\J{4} - \frac{ k_{\ell z} k_{\ell' z}}{\kr^2}\J{5}\Big),\\
    \bs \Theta_{\ell\ell'}^{\downarrow\downarrow} &= \phi_{\ell\ell'}^{\downarrow\downarrow} \Big(\J{1}+ \frac{1}{\kr^2}\J{5}\Big) + \frac{\mu_{\ell} \psi_{\ell\ell' }^{\downarrow\downarrow}}{\mu_{\ell'}k_{\ell}^2} \big(\kr^2\J{2} - \ri k_{\ell z} \J{3} -\ri k_{\ell' z} \J{4} - \frac{k_{\ell z} k_{\ell' z}}{\kr^2} \J{5}\big) ,  \\
    \bs \Theta_{\ell\ell'}^{\downarrow\uparrow} &= \phi_{\ell\ell'}^{\downarrow\uparrow} \Big(\J{1}+ \frac{1}{\kr^2}\J{5}\Big) + \frac{\mu_{\ell} \psi_{\ell\ell'}^{\downarrow\uparrow}}{\mu_{\ell'}k_{\ell}^2}\Big(\kr^2 \J{2} - \ri k_{\ell z} \J{3} +\ri k_{\ell' z} \J{4} + \frac{k_{\ell z} k_{\ell' z}}{\kr^2}  \J{5}\Big) .
\end{split}
\end{equation}
Note that the angular terms in \cref{mat-J} can be rewritten as
\begin{equation}\label{kxky_simply}
    \begin{aligned}
        \frac{k_x}{\kr} &= \frac{e^{\ri \alpha} + e^{-\ri\alpha}}{2}, \quad
        \frac{k_y}{\kr} = \frac{\ri(e^{-\ri \alpha} - e^{\ri\alpha})}{2}, \\
        \frac{k_x^2}{\kr^2} &= \frac{1}{2} + \frac{e^{2\ri\alpha} +e^{-2\ri\alpha}}{4} , \quad
        \frac{k_x k_y}{\kr^2} = \frac{\ri(e^{-2\ri\alpha} - e^{2\ri\alpha})}{4}, \quad
        \frac{k_y^2}{\kr^2} = \frac{1}{2} -\frac{e^{2\ri\alpha} +e^{-2\ri\alpha}}{4}.
    \end{aligned}
\end{equation}
With  
\begin{equation*}
    \begin{split}
        {\bs M}_1=&\begin{bmatrix}
            \frac{1}{2} & 0 & 0\\[7pt]
            0 & \frac{1}{2} & 0\\[7pt]
            0 & 0 & 0
        \end{bmatrix},\quad {\bs M}_2=\begin{bmatrix}
            -\frac{1}{4} & \frac{\ri}{4} & 0\\[7pt]
            \frac{\ri}{4} & \frac{1}{4} & 0\\[7pt]
            0 & 0 & 0
        \end{bmatrix},\quad {\bs M}_3=\begin{bmatrix}
            -\frac{1}{4} & -\frac{\ri}{4} & 0\\[7pt]
            -\frac{\ri}{4} & \frac{1}{4} & 0\\[7pt]
            0 &0 & 0
        \end{bmatrix},\\
        {\bs M}_4=&\begin{bmatrix}
            0 & 0 & \frac{1}{2}\\[7pt]
            0 & 0 & -\frac{\ri}{2}\\[7pt]
            0 & 0 & 0
        \end{bmatrix},\quad {\bs M}_5=\begin{bmatrix}
            0 & 0 & \frac{1}{2}\\[7pt]
            0 & 0 & \frac{\ri}{2}\\[7pt]
            0 & 0 & 0
        \end{bmatrix},\quad {\bs M}_6=\begin{bmatrix}
            0 & 0 & 0\\[7pt]
            0 & 0 & 0\\[7pt]
            0 & 0 & 1
        \end{bmatrix},
    \end{split}
\end{equation*}
the matrices in the expressions \cref{Theta} can be rewritten as 
\begin{equation*}
    \begin{split}
       & \J{1}+ \frac{1}{\kr^2}\J{5}={\bs M}_1+e^{2\ri\alpha}{\bs M}_2+e^{-2\ri\alpha}{\bs M}_3,\quad \J{2}={\bs M}_6,\quad -\frac{\ri\J{3}}{\kr}=e^{\ri\alpha}{\bs M}_4+e^{-\ri\alpha}{\bs M}_5,\\
        &-\frac{\ri}{\kr}\J{4}=e^{\ri\alpha}{\bs M}_4^T+e^{-\ri\alpha}{\bs M}_5^T,\quad  -\frac{1}{\kr^2}\J{5}={\bs M}_1-e^{2\ri\alpha}{\bs M}_2-e^{-2\ri\alpha}{\bs M}_3.
    \end{split}
\end{equation*}
Define $\gamma_{\ell\ell'}=\mu_{\ell}/(\mu_{\ell'}k_{\ell}^2)$ and the following density functions 
	\begin{equation}
		\begin{split}
			\sigma_{\ell\ell'1}^{\uparrow\downarrow}(\kr) &= \frac{\phi_{\ell\ell'}^{\uparrow\downarrow}(\kr)}{k_{\ell'z}}-\gamma_{\ell\ell'}k_{\ell z}\psi_{\ell\ell'}^{\uparrow\downarrow}(\kr),\;\sigma_{\ell\ell'1}^{\downarrow\uparrow} (\kr)= \frac{\phi_{\ell\ell'}^{\downarrow\uparrow}(\kr)}{k_{\ell'z}}-\gamma_{\ell\ell'}k_{\ell z}\psi_{\ell\ell'}^{\downarrow\uparrow}(\kr),\\
			\sigma_{\ell\ell'1}^{\uparrow\uparrow} (\kr)&= \frac{\phi_{\ell\ell'}^{\uparrow\uparrow}(\kr)}{k_{\ell'z}}+\gamma_{\ell\ell'}k_{\ell z}\psi_{\ell\ell'}^{\uparrow\uparrow}(\kr),\; \sigma_{\ell\ell'1}^{\downarrow\downarrow} (\kr)= \frac{\phi_{\ell\ell'}^{\downarrow\downarrow}(\kr)}{k_{\ell'z}}+\gamma_{\ell\ell'}k_{\ell z}\psi_{\ell\ell'}^{\downarrow\downarrow}(\kr),\\
            \sigma_{\ell\ell'3}^{\ast\downarrow} (\kr)&= \frac{\phi_{\ell\ell'}^{\ast\downarrow}(\kr)}{k_{\ell'z}}+\gamma_{\ell\ell'}k_{\ell z}\psi_{\ell\ell'}^{\ast\downarrow}(\kr),\; \sigma_{\ell\ell'3}^{\ast\uparrow}(\kr) = \frac{\phi_{\ell\ell'}^{\ast\uparrow}(\kr)}{k_{\ell'z}}-\gamma_{\ell\ell'}k_{\ell z}\psi_{\ell\ell'}^{\ast\uparrow}(\kr)
		\end{split}
	\end{equation}
	and
	\begin{equation}
		\begin{split}
			\sigma_{\ell\ell'2}^{\uparrow\star}(\kr) &= -\gamma_{\ell\ell'}\frac{k_{\rho}k_{\ell z}}{k_{\ell' z}} \psi_{\ell\ell'}^{\uparrow\star},\quad \sigma_{\ell\ell'2}^{\downarrow\star} (\kr)= \gamma_{\ell\ell'}\frac{k_{\rho}k_{\ell z}}{k_{\ell' z}} \psi_{\ell\ell'}^{\downarrow\star},\\
			\sigma_{\ell\ell'4}^{\ast\downarrow}(\kr)&= \gamma_{\ell\ell'}k_{\rho} \psi_{\ell\ell'}^{\ast\downarrow},\quad\sigma_{\ell\ell'4}^{\ast\uparrow}(\kr) = -\gamma_{\ell\ell'}k_{\rho} \psi_{\ell\ell'}^{\ast\uparrow},\quad
			\sigma_{\ell\ell'5}^{\ast\star} (\kr)= \gamma_{\ell\ell'}\frac{k_{\rho}^2}{k_{\ell' z}} \psi_{\ell\ell'}^{\ast\star},
		\end{split}
	\end{equation}
for all $*,\star=\uparrow,\downarrow$. Substituting into \cref{Theta} and \cref{reactioncomponentspectral} and rearranging the terms lead to the expression
	\begin{equation*}
		\begin{split}
			\widehat{\bd{G}}_{\ell\ell'}^{\ast\star}(k_x, k_y, z, z')=& \frac{\i}{2} {Z}_{\ell\ell'}^{\ast\star}(\sigma_{\ell\ell'1}^{\ast\star}\bs M_1+\sigma_{\ell\ell'3}^{\ast\star}e^{2\ri\alpha}\bs M_2+\sigma_{\ell\ell'3}^{\ast\star}e^{-2\ri\alpha}\bs M_3+\sigma_{\ell\ell'2}^{\ast\star}e^{\ri\alpha}\bs M_4\\
			&+\sigma_{\ell\ell'2}^{\ast\star}e^{-\ri\alpha}\bs M_5+\sigma_{\ell\ell'4}^{\ast\star}e^{\ri\alpha}\bs M_4^{\rm T}+\sigma_{\ell\ell'4}^{\ast\star}e^{-\ri\alpha}\bs M_5^{\rm T}+\sigma_{\ell\ell'5}^{\ast\star}\bs M_6).
		\end{split}
	\end{equation*}
Again, this is a formulation where all the dependence on the polar angle $\alpha$ is found in the $e^{\ri\kappa  \alpha}$ factors, $\kappa=-2,-1,0,1,2$.
Taking inverse Fourier transform gives
	\begin{equation}\label{reactcompphysicaldomain}
		\begin{split}
			\bd{G}_{\ell\ell'}^{\ast\star}{}&{} (\bs{r}, \bs{r}') =\frac{1}{4\pi^2} \iint_{\mathbb{R}^2} \widehat{\bd{G}}_{\ell\ell'}^{\ast\star} e^{\i k_x (x - x') + \i k_y (y - y')} dk_x dk_y\\
			={}&{} \mathcal L_{\ell\ell'0}^{\ast\star}[\sigma_{\ell\ell'1}^{\ast\star}]{\bs M}_1+\mathcal L_{\ell\ell'2}^{\ast\star}[\sigma_{\ell\ell'3}^{\ast\star}]{\bs M}_2+\mathcal L_{\ell\ell',-2}^{\ast\star}[\sigma_{\ell\ell'3}^{\ast\star}]{\bs M}_3+\mathcal L_{\ell\ell'1}^{\ast\star}[\sigma_{\ell\ell'2}^{\ast\star}]{\bs M}_4\\
			+{}&{}\mathcal L_{\ell\ell',-1}^{\ast\star}[\sigma_{\ell\ell'2}^{\ast\star}]{\bs M}_5+\mathcal L_{\ell\ell'1}^{\ast\star}[\sigma_{\ell\ell'4}^{\ast\star}]{\bs M}_4^{\rm T}+\mathcal L_{\ell\ell',-1}^{\ast\star}[\sigma_{\ell\ell'4}^{\ast\star}]{\bs M}_5^{\rm T}+\mathcal L_{\ell\ell'0}^{\ast\star}[\sigma_{\ell\ell'5}^{\ast\star}]{\bs M}_6,
		\end{split}
	\end{equation}
where 
	\begin{equation}\label{generalintegral}
		\mathcal L_{\ell\ell'\kappa}^{\ast\star}[\sigma](\bs r, \bs r')=\frac{\i}{8\pi^2} \iint_{\mathbb{R}^2} {\mathcal E}_{\ell\ell'}^{\ast\star}(\hat\tau_{\ell\ell'}^{\ast\star}(\bs r), \breve\tau_{\ell\ell'}^{\star}(\bs r'))e^{\ri \kappa\alpha}\sigma(k_{\rho}) dk_x dk_y,
	\end{equation}
for $\kappa=-2, -1, 0, 1, 2$, with
	\begin{equation}\label{eq_E}
		{\mathcal E}_{\ell\ell'}^{\ast\star}(\bs r, \bs r')=e^{\i k_x (x - x') + \i k_y (y - y')}e^{s_{\ell\ell'}^{\ast\star} (\ri k_{\ell z}z-\ri k_{\ell' z}z')},
	\end{equation} 
and signs $s_{\ell\ell'}^{\ast\star}$ defined in \cref{eq_signs}. Here, we use the \emph{equivalent polarization coordinates}
\begin{equation}\label{equivalent_imag_for_source}
\breve\tau_{\ell\ell'}^{\downarrow}(\bs r')=\begin{cases}
			\tau_{\ell'}(\bs r'),& \ell\le\ell',\\
			\bs r',& \ell>\ell',
		\end{cases}\quad
		\breve\tau_{\ell\ell'}^{\uparrow}(\bs r')=\begin{cases}
			\bs r',&\ell<\ell',\\
			\tau_{\ell'-1}(\bs r'),&\ell\ge\ell',
		\end{cases}
\end{equation}
for source $\bs r'$ and \emph{effective locations}
\begin{equation}\label{equivalent_imag_for_target}
\begin{split}
&\hat\tau_{\ell\ell'}^{\uparrow\uparrow}(\bs r)=\begin{cases}
\bs r,& \ell<\ell',\\
\tau_{\ell}(\bs r),& \ell\ge\ell',
\end{cases}\quad\hat\tau_{\ell\ell'}^{\uparrow\downarrow}(\bs r)=\begin{cases}
\bs r,& \ell\le\ell',\\
\tau_{\ell}(\bs r),& \ell> \ell',
\end{cases}\\
&\hat\tau_{\ell\ell'}^{\downarrow\uparrow}(\bs r)=\begin{cases}
\tau_{\ell-1}(\bs r),& \ell<\ell',\\
\bs r,& \ell\ge\ell',
\end{cases}\quad \hat\tau_{\ell\ell'}^{\downarrow\downarrow}(\bs r)=\begin{cases}
\tau_{\ell-1}(\bs r),& \ell\le\ell',\\
\bs r,& \ell> \ell',
\end{cases}
\end{split}
\end{equation}
for target $\bs r$ with
\begin{equation}
 \tau_{\ell}(\bs r):=(x, y, \tau_{\ell}(z)),
\end{equation}
where $\tau_{\ell}(\cdot)$ refers to the reflection defined in \cref{z-reflection}.
The formulations \cref{reactcompphysicaldomain,generalintegral} allow us to develop far-field approximations in a uniform framework for all reaction field components $\bd{G}_{\ell\ell'}^{\ast\star}(\bs{r}, \bs{r}') $ regardless of their polar angular dependence.

\subsection{Effective transmission distance of reaction fields}\label{sec_transmission_distance}
\begin{figure}[t]  
		\centering
		\begin{subfigure}{0.38\textwidth}
			\includegraphics[width=\linewidth]{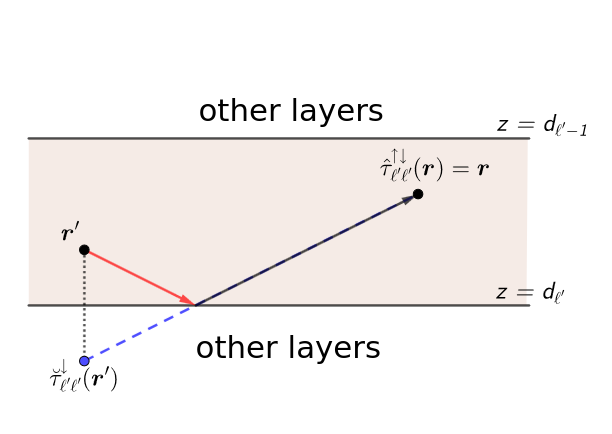}
			\caption{transmission distance for $\widehat{\bs G}_{\ell'\ell'}^{\uparrow\downarrow}$}\label{fig_targetoversource_eq_a}
		\end{subfigure}
		\quad
		\begin{subfigure}{0.41\textwidth}
			\includegraphics[width=\linewidth]{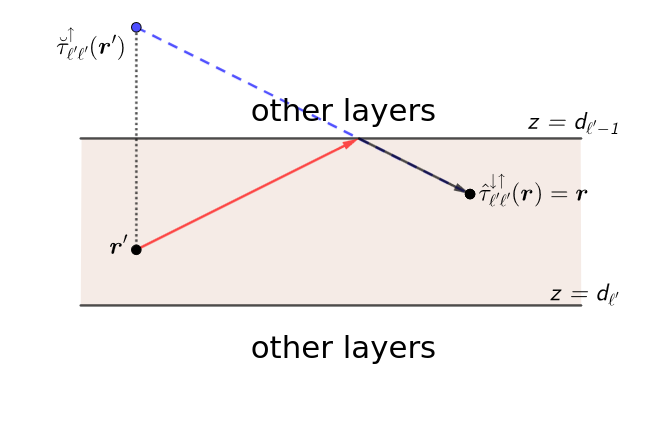}
			\caption{transmission distance for $\widehat{\bs G}_{\ell'\ell'}^{\downarrow\uparrow}$}\label{fig_targetoversource_eq_b}
		\end{subfigure}
        \begin{subfigure}{0.4\textwidth}
			\includegraphics[width=\linewidth]{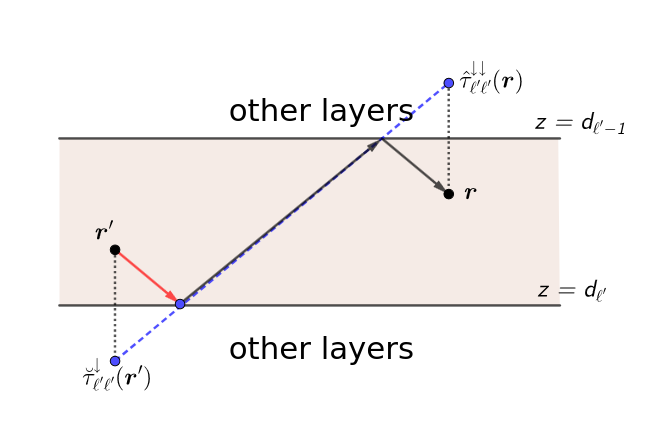}
			\caption{transmission distance for $\widehat{\bs G}_{\ell'\ell'}^{\downarrow\downarrow}$}\label{fig_targetoversource_eq_c}
		\end{subfigure}
		\quad
		\begin{subfigure}{0.4\textwidth}
			\includegraphics[width=\linewidth]{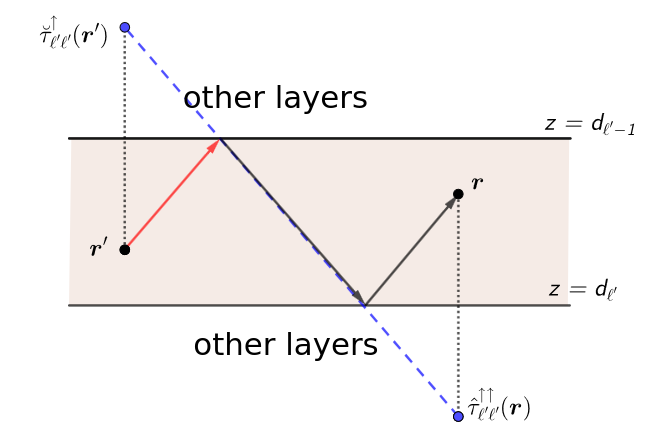}
			\caption{transmission distance for $\widehat{\bs G}_{\ell'\ell'}^{\uparrow\uparrow}$}\label{fig_targetoversource_eq_d}
		\end{subfigure}
		\caption{Equivalent polarization source coordinates and effective target locations in the case of $\ell=\ell'$.}
		\label{targetoversource_eq}
\end{figure}

\begin{figure}[t]  
\centering
\begin{subfigure}{0.38\textwidth}
\includegraphics[width=\linewidth]{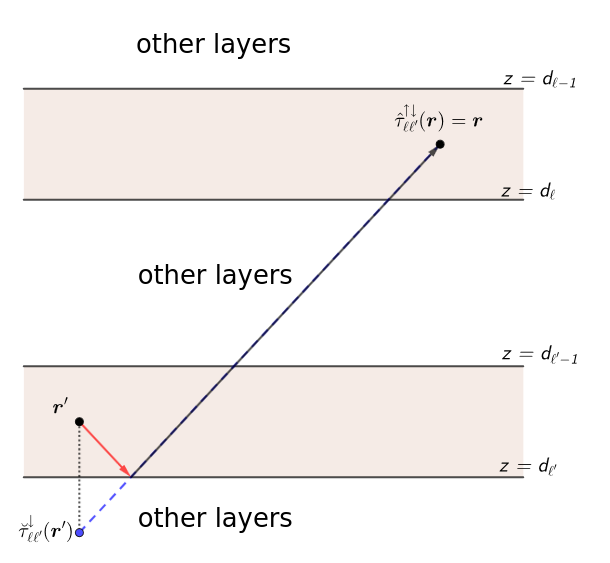}
\caption{transmission distance for $\widehat{G}_{\ell\ell'}^{\uparrow\downarrow}$}
\end{subfigure}\qquad
\begin{subfigure}{0.4\textwidth}
\includegraphics[width=\linewidth]{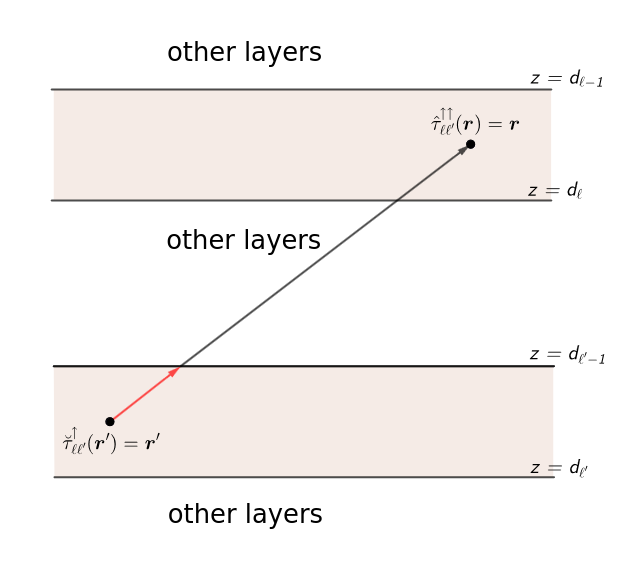}
\caption{transmission distance for $\widehat{G}_{\ell\ell'}^{\uparrow\uparrow}$}
\end{subfigure}
\begin{subfigure}{0.43\textwidth}
\includegraphics[width=\linewidth]{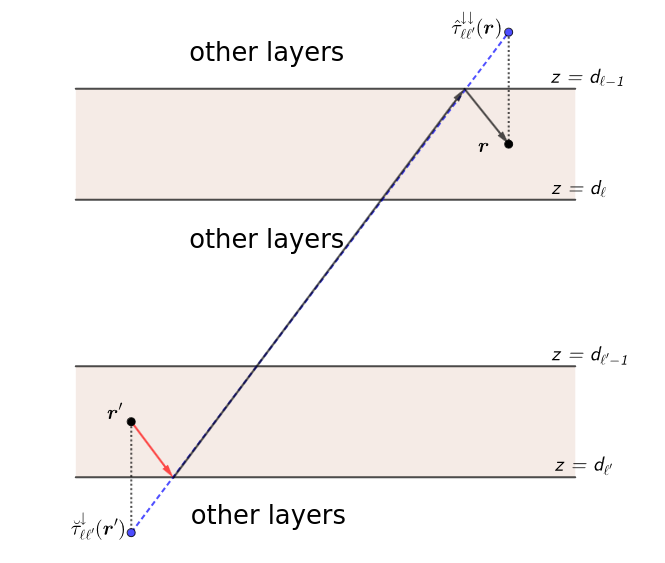}
\caption{transmission distance for $\widehat{G}_{\ell\ell'}^{\downarrow\downarrow}$}
\end{subfigure}
\quad
\begin{subfigure}{0.42\textwidth}
\includegraphics[width=\linewidth]{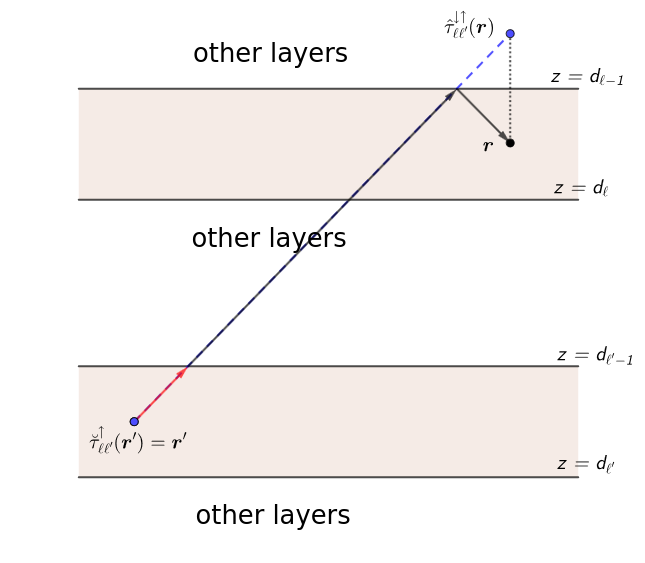}
\caption{transmission distance for $\widehat{G}_{\ell\ell'}^{\downarrow\uparrow}$}
\end{subfigure}
\caption{Equivalent polarization source coordinates and effective target locations in the case of $\ell<\ell'$.}
\label{fig_targetoversource}
\end{figure}

\begin{figure}[ht!]  
\centering
\begin{subfigure}{0.4\textwidth}
\includegraphics[width=\linewidth]{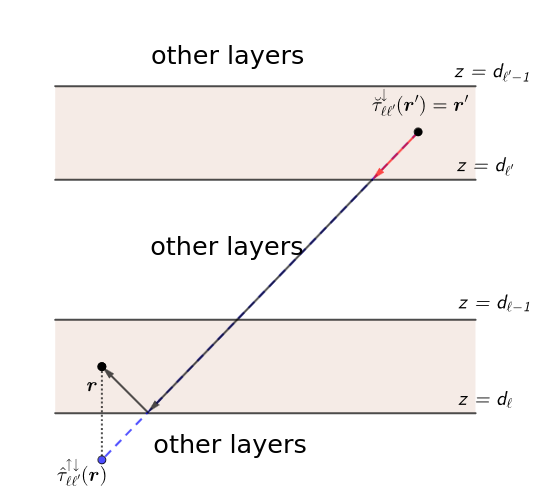}
\caption{transmission distance for $\widehat{G}_{\ell\ell'}^{\uparrow\downarrow}$}
\end{subfigure}
\quad
\begin{subfigure}{0.38\textwidth}
\includegraphics[width=\linewidth]{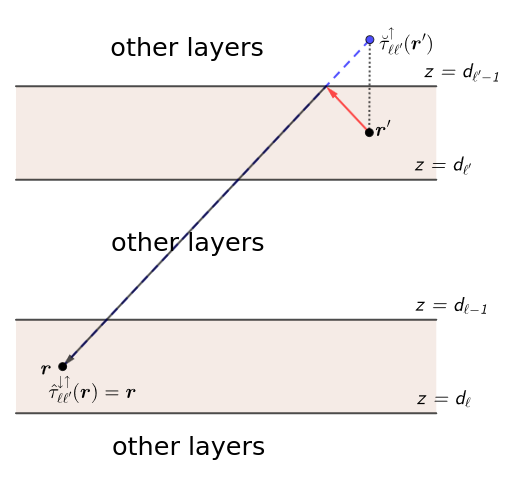}
\caption{transmission distance for $\widehat{G}_{\ell\ell'}^{\downarrow\uparrow}$}
\end{subfigure}
\begin{subfigure}{0.4\textwidth}
\includegraphics[width=\linewidth]{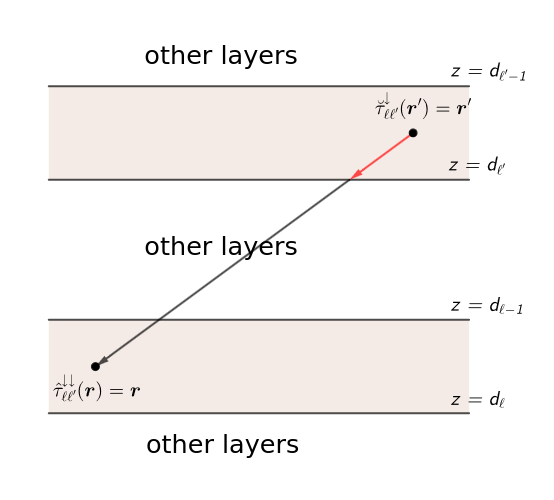}
\caption{transmission distance for $\widehat{G}_{\ell\ell'}^{\downarrow\downarrow}$}
\end{subfigure}
\quad
\begin{subfigure}{0.38\textwidth}
\includegraphics[width=\linewidth]{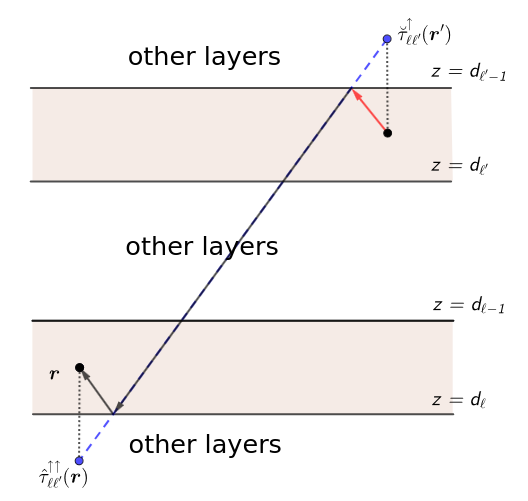}
\caption{transmission distance for $\widehat{G}_{\ell\ell'}^{\uparrow\uparrow}$}
\end{subfigure}
\caption{Equivalent polarization source coordinates and effective target locations in the case of $\ell>\ell'$.}
\label{fig_sourceovertarget}
\end{figure}

Besides the equivalent polarization coordinates proposed in our previous works \cite{bo2019hfmm,zhang2020}, we introduce a new concept of effective location for the target particles, in order to account for the actual transmission distance of the reflected waves in layered media.
The upward and downward waves generated by the source at $\bs r'$ transmit to the target at $\bs r$ via different paths, see \cref{targetoversource_eq,fig_targetoversource,fig_sourceovertarget}, and generally induce upward and downward reaction fields. That's the physical background of the decomposition of four reaction field components.
The reaction field decomposition \cref{componets} of the dyadic Green's functions consists of the free-space interaction $\widehat{\bd{G}}_{\bs E}^{f}$ provided $\ell = \ell'$, as well as the reaction field.
When categorized by the upward/downward field propagation directions, the reaction field is decomposed into (up to) four terms in \cref{componets}, each $\widehat{\bd{G}}_{\ell\ell'}^{\ast\star}$ representing one type with upward ($\uparrow$), downward ($\downarrow$) or both directions, with the first symbol $\ast$ indicating the direction of wave arriving at the target, and the second symbol  $\star$  indicating the direction of the wave leaving the source.

When the target and the source come from the same layer, i.e. $\ell = \ell'$, waves of the reaction field must have at least one reflection on interfaces due to the subtraction of the free-space part, see \cref{targetoversource_eq} for an illustration. For instance, in \cref{fig_targetoversource_eq_c}, the reaction field component $\widehat{\bd{G}}_{\ell'\ell'}^{\downarrow\downarrow}$ is interpreted as the superposition of waves that are downward at $\bs r'$ and downward at $\bs r$, including the wave marked by the solid line with two reflections in the figure, as well as any other contributions that may have experienced more reflections and transmissions on the interfaces. The minimal vertical transmission distance of these waves is given by the distance between the effective target location $\hat\tau_{\ell'\ell'}^{\downarrow\downarrow}(\bs r) = \tau_{\ell'-1}(\bs r)$ and the equivalent polarization source $\breve\tau_{\ell'\ell'}^{\downarrow}(\bs r') = \tau_{\ell'}(\bs r')$, as shown by the dashed line.
Indeed, in the exponent of \cref{zexponential},
\begin{equation*}
Z^{\downarrow\downarrow}_{\ell'\ell'}(\kr, z, z') \sim e^{-\kr(\tau_{\ell'-1}(z) - \tau_{\ell'}(z'))}, \quad \kr \to +\infty.
\end{equation*}

The cases which the target and the source are located in different layers are illustrated by \cref{fig_targetoversource} for $\ell < \ell'$, and \cref{fig_sourceovertarget} for $\ell > \ell'$, respectively.
The only difference with the previous case is that one of the reaction field components has minimal vertical transmission distance given by $|z - z'|$. This component must exist because it is \emph{not} equivalent to the free-space interaction due to the transmission through multiple layers.

In later discussion of the FMM implementation, the field transmission distance
\begin{equation}\label{eq_transdict}
    d_{\ell\ell'}^{\ast\star}(\bs{r}, \bs{r}') = \left| \hat\tau_{\ell\ell'}^{\ast\star}(\bs r) - \breve\tau_{\ell\ell'}^{\star}(\bs r') \right|
\end{equation}
will be used as the criterion of far-field expansions.
Note that the field transmission distance is not shorter than that based on equivalent polarization source alone in our previous works \cite{bo2019hfmm,zhang2020}  for handling Helmholtz equation in layered media, suggesting that wave sources in layered media are even more separated than previously thought.

\subsection{FMM in layered media}
Let $\mathcal{P}_{\ell} = {(\bs q_{\ell j}, \bs r_{\ell j})}_{j=1}^{N_{\ell}}$ for $\ell = 0, 1, \dots, L$ denote $L+1$ groups of wave sources, where each group is located in the $\ell$-th layer of a multilayered medium consisting of $L+1$ layers. The interactions between all $N_{tot} = N_0 + N_1 + \cdot \cdot \cdot + N_L$ particles are given by the summations
\begin{equation}\label{layeredmediafields}
	\bs \Phi_{\ell'}(\bs r_{\ell' i})=\begin{bmatrix}
		\bs \Phi_{\ell'x}(\bs r_{\ell' i})\\
		\bs \Phi_{\ell'y}(\bs r_{\ell' i})\\
		\bs \Phi_{\ell'z}(\bs r_{\ell' i})
		\end{bmatrix}=\bs\Phi^{f}_{\ell'}(\bs r_{\ell' i})+\sum_{\ell=0}^{L}\sum_{\ast,\star=\uparrow,\downarrow}\bs\Phi_{\ell\ell'}^{\ast\star}(\bs r_{\ell' i}),
\end{equation}
for $\ell=0,1,\cdots,L$ and $i=1,2,\cdots,N_{\ell}$, where 
\begin{equation}\label{allcomponents}
	\bs\Phi^{f}_{\ell'}(\bs r_{\ell' i})=\sum_{j=1,j\not=i}^{N_{\ell'}}\bs{G}_{\bs E}^f(\bs r_{\ell' j},\bs r_{\ell' i})\bs q_{\ell' j},\quad \bs\Phi_{\ell\ell'}^{\ast\star}(\bs r_{\ell' i})=\sum_{j=1}^{N_{\ell}}\bs{G}_{\ell\ell'}^{\ast\star}(\bs r_{\ell j},\bs r_{\ell' i})\bs q_{\ell j}.
\end{equation}
Here, we put the summation on the first coordinates $\bs r_{\ell j}$ for the sake of future application in accelerating the integral methods (e.g., method of moments) for solving electromagnetic problems. It is equivalent to the summation on the second coordinates $\bs r_{\ell' i}$ due to the symmetry \cite{tai1971dyadic} of the dyadic Green's function. 

Like in the previous works \cite{bo2020tefmm,bo2019hfmm}, we separately implement the FMM for the free-space part using classic approaches, and for each reaction field component marked by the quadruple $(\ell, \ell', \ast, \star)$, then sum up the results. In the low and medium frequency regime, the overall computational cost is $\mathcal{O}(LN_{tot}\log N_{tot})$.
Without loss of generality, we assume the same number $N$ of targets in layer $\ell$ and sources in layer $\ell'$, respectively. Consider the reaction field component
\begin{equation}\label{reactionfields}
	\bs \Phi(\bs r_i)=\bs\Phi_{\ell\ell'}^{\ast\star}(\bs r_i)=
	\begin{bmatrix}
	\Phi_x(\bs r_i)\\
	\Phi_y(\bs r_i)\\
	\Phi_z(\bs r_i)
	\end{bmatrix}:=\sum\limits_{j=1}^{N}\bd{G}_{\ell\ell'}^{\ast\star}(\bs{r}_j, \bs{r}_i)\begin{bmatrix}
		q_{j}^x\\
		q_{j}^y\\
		q_{j}^z
	\end{bmatrix}
\end{equation}
Henceforth, the dependence on $\ell$, $\ell'$, and propagating direction notations $\ast,\star$ will be omitted for simplicity.

\begin{figure}[!ht]  
    \centering
    \includegraphics[width=0.9\linewidth]{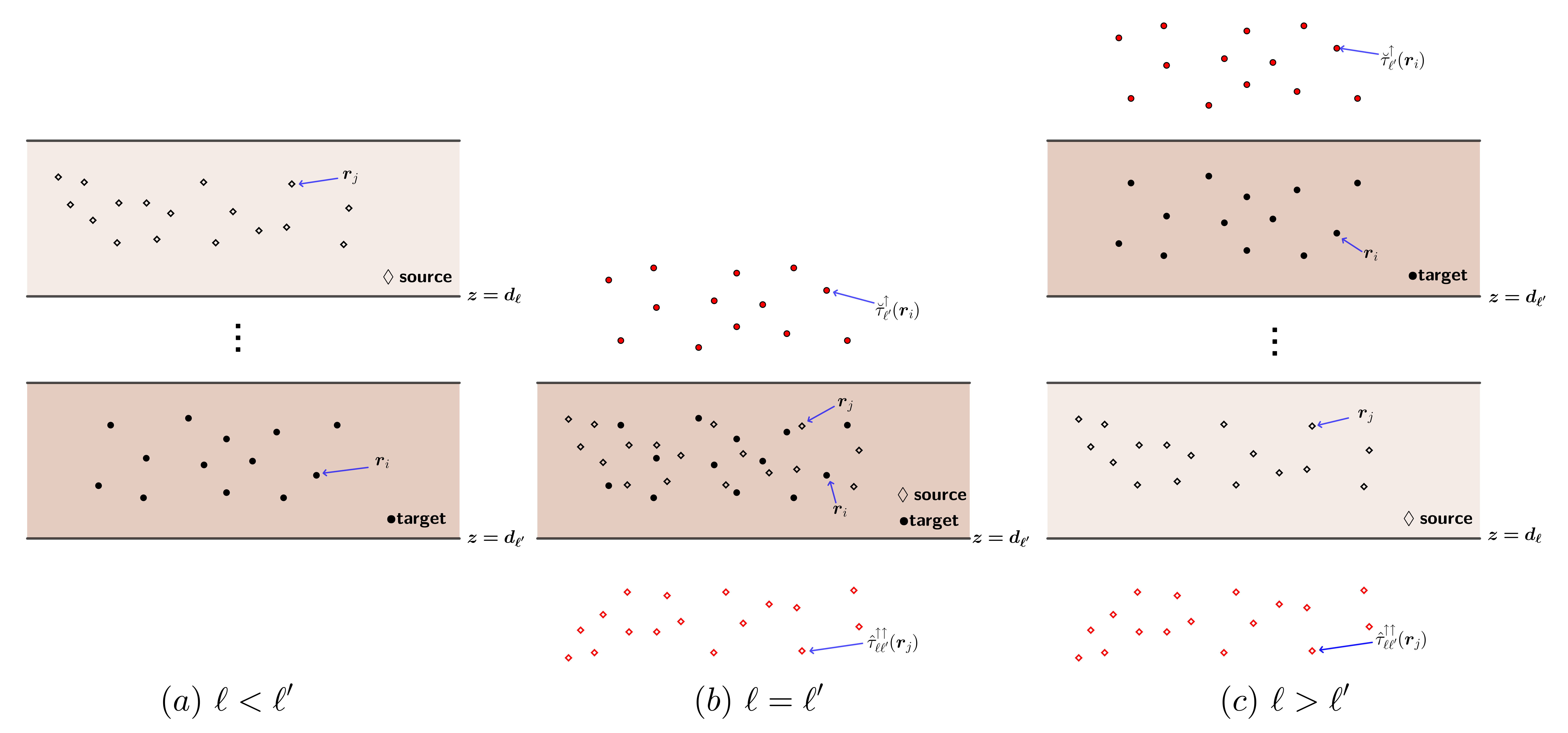}
    \caption{The effective locations $\{\breve\tau_{\ell'}^{\uparrow}(\bs r_i)\}$ and the equivalent polarization coordinates $\{\hat\tau_{\ell\ell'}^{\uparrow\uparrow}(\bs r_j)\}$.}
    \label{equivalent_particles_withuparrow}
\end{figure}
\begin{figure}[!ht]  
    \centering
    \includegraphics[width=0.9\linewidth]{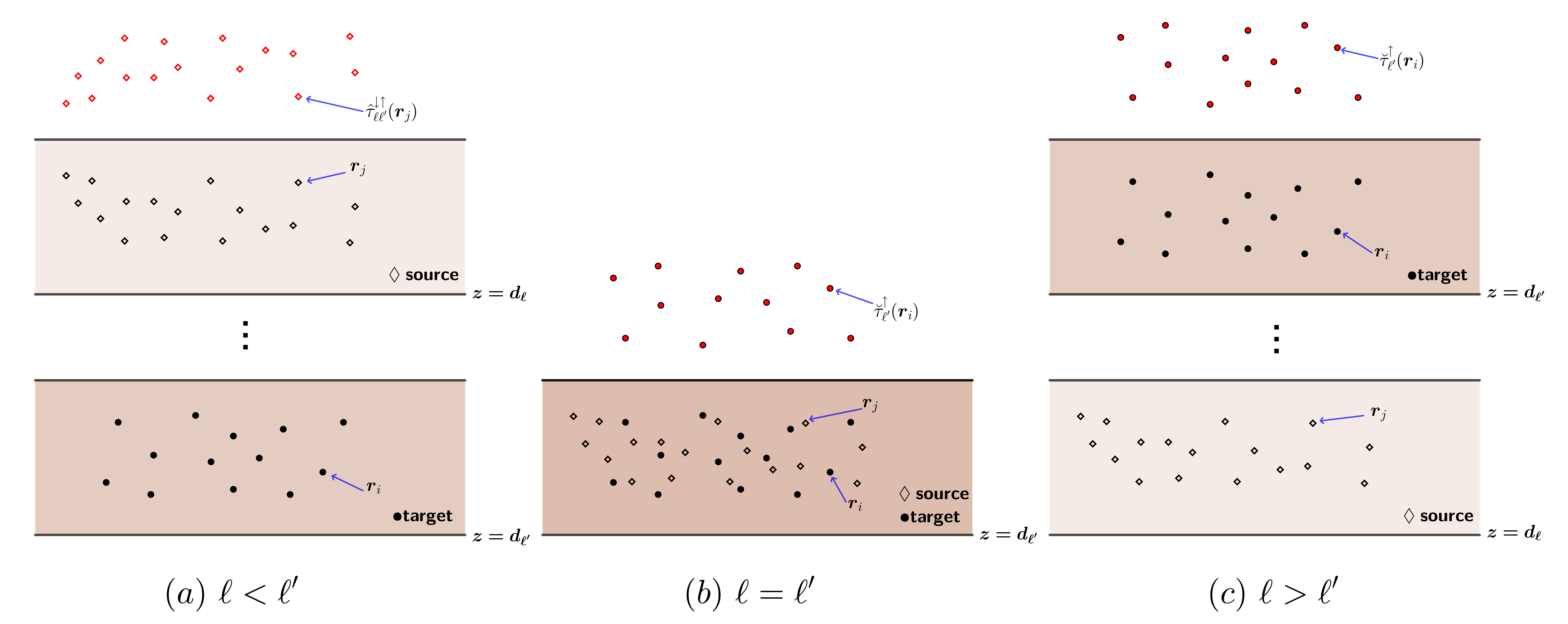}
    \caption{The effective locations $\{\breve\tau_{\ell'}^{\uparrow}(\bs r_{i})\}$ and equivalent polarization coordinates $\{\hat\tau_{\ell\ell'}^{\downarrow\uparrow}(\bs r_{j})\}$.}
    \label{equivalent_particles_withdownarrowuparrow}
\end{figure}

In the integral representation \cref{generalintegral}, we have applied the equivalent polarization coordinates $\hat\tau_{\ell\ell'}^{\ast\star}(\bs r_j)$ and effective locations $\breve\tau_{\ell'}^{\star}(\bs r_i)$, which are transforms on the target layer and the source layer, respectively, see \cref{equivalent_particles_withuparrow} and \cref{equivalent_particles_withdownarrowuparrow}, which show the upward direction of the wave leaving the source. The transformed layers are always separated by an interface plane.

To simplify the notations, we omit the subscripts $\ell$, $\ell'$ and denote by
\begin{equation*}
\begin{split}
    k&=k_{\ell},\enspace k'=k_{\ell'},\enspace k_z=k_{\ell z},\enspace k_z'=k_{\ell' z},\enspace \overline{\bs r}_{j}=\hat\tau_{\ell\ell'}^{\ast\star}(\bs r_j), \enspace \overline{\bs r}_{i}^{\prime}=\breve\tau_{\ell'}^{\star}(\bs r_i), \enspace \sigma_{m}=\sigma_{\ell\ell'm}^{\ast\star},\\
    \mathcal{Z}&(k_{\rho},z,z')=e^{s_{\ell\ell'}^{\ast\star} (\ri k_{z}z-\ri k_{z}'z')},\quad {\mathcal E}(\bs r, \bs r')=e^{\ri k_x (x - x') + \ri k_y (y - y')}\mathcal{Z}(k_{\rho},z,z').
\end{split}
\end{equation*}
Define integrals
\begin{equation}\label{general_integral_I}
	\mathcal{I}_{nm\kappa}^{\nu \mu}(\bs r,\bs r',\sigma)=\dfrac{c_{nm}^{\nu\mu}}{2\pi}\iint_{\mathbb{R}^2}{\mathcal E}(\bs r, \bs r')e^{\ri(m-\mu+\kappa)\alpha}\widehat{P}_{n}^{m}\Big(\frac{k_{z}}{k}\Big)\widehat{P}_{\nu}^{\mu}\Big(\frac{k_{z}'}{k'}\Big)\sigma(k_{\rho}) dk_x dk_y,
\end{equation}
where $c_{nm}^{\nu\mu}=(s_{\ell\ell'}^{\ast\star})^{n+\nu+m+\mu}\i^{n+\nu+1}$, and $\widehat P_n^m(x)$ is the analytic extension of the normalized associated Legendre function.
In particular,
\begin{equation}\label{general_integral_I_00}
	\mathcal{I}_{00\kappa}^{00}(\bs r,\bs r',\sigma)=\frac{\i}{8\pi^2} \iint_{\mathbb{R}^2}{\mathcal E}_{\ell\ell'}^{\ast\star}(\bs r, \bs r')e^{\ri\kappa\alpha}\sigma(k_{\rho}) dk_x dk_y=\mathcal L_{\ell\ell'\kappa}^{\ast\star}[\sigma](\bs r, \bs r').
\end{equation}
By \cref{reactcompphysicaldomain,allcomponents,reactionfields}, we obtain
\begin{equation}\label{reactcompphysicaldomainmat}
    \begin{split}
		\bs \Phi(\bs r_i)={}&{}{\bs M}_1\sum\limits_{j=1}^N\bs q_j\mathcal{I}_{000}^{00}(\overline{\bs r}_j, \overline{\bs r}_i',\sigma_{1})+{\bs M}_6\sum\limits_{j=1}^N\bs q_j\mathcal{I}_{000}^{00}(\overline{\bs r}_j, \overline{\bs r}_i',\sigma_{5})\\
		{}&{}+{\bs M}_5\sum\limits_{j=1}^N\bs q_j\mathcal{I}_{00,-1}^{00}(\overline{\bs r}_j, \overline{\bs r}_i',\sigma_{2})+{\bs M}_5^{\rm T}\sum\limits_{j=1}^N\bs q_j\mathcal{I}_{00,-1}^{00}(\overline{\bs r}_j, \overline{\bs r}_i',\sigma_{4})\\
		{}&{}+{\bs M}_4\sum\limits_{j=1}^N\bs q_j\mathcal{I}_{001}^{00}(\overline{\bs r}_j, \overline{\bs r}_i',\sigma_{2})+{\bs M}_4^{\rm T}\sum\limits_{j=1}^N\bs q_j\mathcal{I}_{001}^{00}(\overline{\bs r}_j, \overline{\bs r}_i',\sigma_{4})\\
		{}&{}+{\bs M}_3\sum\limits_{j=1}^N\bs q_j\mathcal{I}_{00,-2}^{00}(\overline{\bs r}_j, \overline{\bs r}_i',\sigma_{3})+{\bs M}_2\sum\limits_{j=1}^N\bs q_j\mathcal{I}_{002}^{00}(\overline{\bs r}_j, \overline{\bs r}_i',\sigma_{3}).
    \end{split}
\end{equation}
In \cite{bo2019hfmm}, the far-field expansions of \cref{general_integral_I_00} with $\kappa=0$ have been proposed and numerically verified for exponential convergence. Here, we extend the results for $\kappa = \pm 1, \pm 2$.
By the extended Funk--Hecke formula \cite[Proposition 6]{bo2019hfmm}, we have
\begin{equation*}\label{me_pre}
    \begin{aligned}
		{\mathcal E}(\bs r, \bs r')&={\mathcal E}(\bs r_c, \bs r')e^{\ri k_x (x - x_c) + \ri k_y (y - y_c) +s_{\ell\ell'}^{\ast\star} \ri k_{z}(z-z_c)} \\
		&= {\mathcal E}(\bs r_c, \bs r') \sum\limits_{n=0}^{\infty}\sum\limits_{m=-n}^n4\pi j_n(k r_s)\overline{Y_n^m({\theta}_s,\varphi_s)}\ri^n(s_{\ell\ell'}^{\ast\star})^{n+m}\widehat{P}_n^m\Big(\frac{k_{z}}{k}\Big)e^{\ri m\alpha},
    \end{aligned}
\end{equation*}
and
\begin{equation*}\label{le_pre}
    \begin{aligned}
		{\mathcal E}(\bs r, \bs r')&={\mathcal E}(\bs r, \bs r_c)e^{\ri k_x (x_c - x') + \ri k_y (y_c - y') +s_{\ell\ell'}^{\ast\star} \ri k_{z}'(z_c-z')} \\
	&= {\mathcal E}(\bs r, \bs r_c) \sum\limits_{n=0}^{\infty}\sum\limits_{m=-n}^n4\pi j_n(k' r_s'){Y_n^m({\theta}_s',\varphi_s')}(-\ri)^n(s_{\ell\ell'}^{\ast\star})^{n+m}\widehat{P}_n^m\Big(\frac{k_{z}'}{k'}\Big)e^{-\ri m\alpha},
    \end{aligned}
\end{equation*}
where $({r}_s,{\theta}_s,\varphi_s)$, $({r}_s',{\theta}_s',\varphi_s')$ are the spherical coordinates of $\bs r-\bs r_c$ and $\bs r'-\bs r_c$, respectively.
Substituting the pair of expansions into the integral \cref{general_integral_I_00} and exchanging the order of integration and summation (over $n, \nu$), we obtain the \emph{prototype} expansions
\begin{equation}\label{I_addition_1}
	\mathcal{I}_{00\kappa}^{\nu \mu}(\bs r,\bs r',\sigma)=\sum\limits_{n=0}^{\infty}\sum\limits_{m=-n}^n\sqrt{4\pi}j_n(k r_s)\overline{Y_n^m({\theta}_s,\varphi_s)}\mathcal{I}_{nm\kappa}^{\nu \mu}(\bs r_c,\bs r',\sigma),
\end{equation}
and
\begin{equation}\label{I_addition_2}
	\mathcal{I}_{nm\kappa}^{00}(\bs r,\bs r',\sigma)=\sum\limits_{\nu=0}^{\infty}\sum\limits_{\mu=-\nu}^{\nu}(-1)^{\nu}\sqrt{4\pi} j_{\nu}(k' r_s'){Y_{\nu}^{\mu}({\theta}_s',\varphi_s')}\mathcal{I}_{nm\kappa}^{\nu \mu}(\bs r,\bs r_c,\sigma),
\end{equation}
respectively.
For a general component in \cref{reactcompphysicaldomainmat}, expansion \cref{I_addition_1} implies the ME 
\begin{equation}\label{integralME}
    \sum_{j=1}^{N} q_j^{v}\mathcal{I}_{00\kappa}^{00}(\overline{\bs r}_j, \overline{\bs r}_i',\sigma) 
	=\dfrac{1}{\sqrt{4\pi}}\sum_{n=0}^{\infty}\sum\limits_{m=-n}^n  M_{nm}^v(\bs r_c^s)\mathcal{I}_{nm\kappa}^{00}( \bs{r}_c^s,\overline{\bs r}_i',\sigma),
\end{equation}
for all $v=x, y, z, \kappa=0, \pm 1, \pm 2$ at source center $\bs r_c^s$, where 
\begin{equation}\label{mecoeff}
	M_{ nm}^v(\bs r_c^s) = \sum_{j=1}^{N} 4\pi q_j^{v} j_n(k {\hat r}_{j}) \overline{Y_n^m({\hat\theta}_{j}, {\hat\phi}_{j})},\quad v=x, y, z,
\end{equation}
and $(\hat r_j, \hat\theta_j,\hat\phi_j)$ and $(\hat r_i', \hat\theta_i',\hat\phi_i')$ are the spherical coordinates of $\overline{\bs r}_j-\bs r_c^s$ and $\overline{\bs r}_i'-\bs r_c^s$, respectively.
Similarly, \cref{I_addition_2} implies the LE
\begin{equation}\label{integralLE}
	\sum_{j=1}^{N} q_j^{v}\mathcal{I}_{00\kappa}^{00}(\overline{\bs r}_j, \overline{\bs r}_i',\sigma) 
	=\sum_{n=0}^{\infty}\sum\limits_{m=-n}^nL_{nm\kappa}^v( \bs r_c^t,)j_n(k{'}\tilde r_i')Y_n^m(\tilde\theta_i',\tilde\varphi_i'),
\end{equation}
for all $v=x, y, z, \kappa=0, \pm 1, \pm 2$ at target center $\bs r_c^t$, where 
\begin{equation}\label{LEcoefficient}
    L_{nm\kappa }^v(\bs r_c^t) = (-1)^n\sqrt{4\pi}\sum_{j=1}^{N} q_j^{v} \mathcal{I}_{00\kappa}^{nm}(\overline{\bs r}_j, \bs r_c^t,\sigma),
\end{equation}
and $(\tilde r_j,\tilde \theta_j, \tilde\phi_j)$ and $(\tilde r_i',\tilde\theta_i', \tilde\phi_i')$ are the spherical coordinates of $\overline{\bs r}_j-\bs r_c^t$ and $\overline{\bs r}_i'-\bs r_c^t$, respectively.
Applying \cref{I_addition_1} to the integrals above, we obtain the M2L translation
\begin{equation}\label{metoleimage1}
\begin{split}
L_{nm\kappa }^{v}(\bs r_c^t)&=(-1)^n4\pi\sum_{j=1}^{N} q_j^{v} \sum\limits_{n'=0}^{\infty}\sum\limits_{m'=-n'}^{n'}j_n(k \hat r_j)\overline{Y_n^m(\hat{\theta}_j,\hat{\varphi}_j)}\mathcal{I}_{n'm'\kappa}^{nm}({\bs r}_c^s, \bs r_c^t,\sigma)\\
&=\sum\limits_{n'=0}^{\infty}\sum\limits_{m'=-n'}^{n'}(-1)^n\mathcal{I}_{n'm'\kappa}^{nm}(\bs r_c^s, \bs r_c^t,\sigma)M_{n'm'}^{v}(\bs r_c^s)
\end{split}
\end{equation}
from the ME at $\bs r_c^s$ to the LE at $\bs r_c^t$, for all $v=x, y, z, \kappa=0, \pm 1, \pm 2$.

Applying the ME formula \cref{integralME} to all the integrals in \cref{reactcompphysicaldomainmat}, we obtain the ME of the entire reaction field component
\begin{equation}\label{maxwell_ME_layered}
    \begin{split}
		\bs \Phi(\bs r_i)
		&=\dfrac{1}{\sqrt{4\pi}}\sum_{n=0}^{\infty}\sum\limits_{m=-n}^n\bs F_{nm}(\bs r_c^s, \overline{\bs r}_i')\bs m_{nm},\quad \bs m_{nm}=\begin{bmatrix}
			M_{nm}^x\\
			M_{nm}^{y}\\
			M_{nm}^{z}
		\end{bmatrix},
    \end{split}
\end{equation}
where $M_{nm}^{x}, M_{nm}^y, M_{nm}^z$ are ME coefficients given in \cref{mecoeff}, and
\begin{equation}\label{M2T_F_nm_mg0}
	\begin{aligned}
		\bs F_{nm}(\bs r_c^s, \overline{\bs r}_i')
		={}&{}\mathcal{I}_{nm0}^{00}(\bs r_c^s, \overline{\bs r}_i',\sigma_{1}){\bs M}_1+\mathcal{I}_{nm0}^{00}(\bs r_c^s, \overline{\bs r}_i',\sigma_{5}){\bs M}_6\\
		{}&{} +\mathcal{I}_{nm,-1}^{00}( \bs r_c^s, \overline{\bs r}_i',\sigma_{2}){\bs M}_5+\mathcal{I}_{nm,-1}^{00}(\bs r_c^s, \overline{\bs r}_i',\sigma_{4}){\bs M}_5^{\rm T}\\
		{}&{} +\mathcal{I}_{nm1}^{00}(\bs r_c^s, \overline{\bs r}_i',\sigma_{2}){\bs M}_4+\mathcal{I}_{nm1}^{00}(\bs r_c^s, \overline{\bs r}_i',\sigma_{4}){\bs M}_4^{\rm T}\\
		{}&{} +\mathcal{I}_{nm,-2}^{00}(\bs r_c^s, \overline{\bs r}_i',\sigma_{3}){\bs M}_3+\mathcal{I}_{nm2}^{00}(\bs r_c^s, \overline{\bs r}_i',\sigma_{3}){\bs M}_2
	\end{aligned}
\end{equation}
are the ME basis functions.
Similarly, applying the LE \cref{integralLE} to all the integrals in \cref{reactcompphysicaldomainmat}, we obtain the LE of the entire reaction field component
\begin{equation}\label{maxwell_LE_layered}
	\bs \Phi(\bs r_i)=\sum_{n=0}^{\infty}\sum\limits_{m=-n}^n(-1)^n\sqrt{4\pi}\bs l_{nm}(\bs r_c^t)j_n(k_{\ell'}\tilde r_{i}')Y_n^m(\tilde\theta_i',\tilde\phi_i'),
\end{equation}
where the LE coefficient vectors
\begin{equation}\label{S2L_L_nm_mg0}
    \begin{aligned}
        \bs l_{nm}(\bs r_c^t)
            ={}&{}\sum\limits_{j=1}^N\mathcal{I}_{000}^{nm}(\overline{\bs r}_j, \bs r_c^t, \sigma_{1}){\bs M}_1\bs q_j+\sum\limits_{j=1}^N\mathcal{I}_{000}^{nm}(\overline{\bs r}_j, \bs r_c^t,\sigma_{5}){\bs M}_6\bs q_j\\
		{}&{} +\sum\limits_{j=1}^N\mathcal{I}_{00,-1}^{nm}( \overline{\bs r}_j, \bs r_c^t,\sigma_{2}){\bs M}_5\bs q_j+\sum\limits_{j=1}^N\mathcal{I}_{00,-1}^{nm}(\overline{\bs r}_j, \bs r_c^t,\sigma_{4}){\bs M}_5^{\rm T}\bs q_j\\
		{}&{} +\sum\limits_{j=1}^N\mathcal{I}_{001}^{nm}(\overline{\bs r}_j, \bs r_c^t,\sigma_{2}){\bs M}_4\bs q_j+\sum\limits_{j=1}^N\mathcal{I}_{001}^{nm}( \overline{\bs r}_j,\bs r_c^t,\sigma_{4}){\bs M}_4^{\rm T}\bs q_j\\
		{}&{} +\sum\limits_{j=1}^N\mathcal{I}_{00,-2}^{nm}(\overline{\bs r}_j, \bs r_c^t,\sigma_{3}){\bs M}_3\bs q_j+\sum\limits_{j=1}^N\mathcal{I}_{002}^{nm}(\overline{\bs r}_j, \bs r_c^t,\sigma_{3}){\bs M}_2\bs q_j.
    \end{aligned}
\end{equation}

Finally, we emphasize that the formulations of the ME coefficients in \cref{mecoeff} and the LE basis functions in \cref{maxwell_LE_layered} are precisely the same as those employed in the scalar case for Helmholtz equation \cite{bo2019hfmm}. Accordingly, the multipole-to-multipole (M2M) and local-to-local (L2L) translation operators remain unchanged, and are given by \cref{maxwell_M2M_freespace,maxwell_L2L_freespace}, respectively. Moreover, the M2L translation \cref{metoleimage1} implies
\begin{equation}\label{M2L_scale_with_Tnmu}
	\bs l_{nm}(\bs r_c^t)=(-1)^n\sum\limits_{n'=0}^{\infty}\sum_{m'=-n'}^{n'}\bs T_{nm}^{n'm'}(\bs r_c^s, \bs r_c^t)\bs m_{n'm'}(\bs r_c^s),
\end{equation}
where the translation tensors
\begin{equation}\label{Tnm_nmu_M2L}
	\begin{split}
		\bs T_{nm}^{n'm'}(\bs r_c^s, \bs r_c^t)={}&{}\mathcal{I}_{n'm'0}^{nm}(\bs r_c^s, \bs r_c^t,\sigma_{1}){\bs M}_1+\mathcal{I}_{n'm'0}^{nm}(\bs r_c^s, \bs r_c^t,\sigma_{5}){\bs M}_6\\
		{}&{} +\mathcal{I}_{n'm',-1}^{nm}( \bs r_c^s, \bs r_c^t,\sigma_{2}){\bs M}_5+\mathcal{I}_{n'm',-1}^{nm}(\bs r_c^s, \bs r_c^t,\sigma_{4}){\bs M}_5^{\rm T}\\
		{}&{} +\mathcal{I}_{n'm'1}^{nm}(\bs r_c^s, \bs r_c^t,\sigma_{2}){\bs M}_4+\mathcal{I}_{n'm'1}^{nm}(\bs r_c^s, \bs r_c^t,\sigma_{4}){\bs M}_4^{\rm T}\\
		{}&{} +\mathcal{I}_{n'm',-2}^{nm}(\bs r_c^s, \bs r_c^t,\sigma_{3}){\bs M}_3+\mathcal{I}_{n'm'2}^{nm}(\bs r_c^s, \bs r_c^t,\sigma_{3}){\bs M}_2.
	\end{split}
\end{equation}

Applying the ME, LE, and M2L derived above and the M2M, L2L employed in the scalar Helmholtz equation in the framework proposed in \cite{bo2019hfmm} implements the FMM for $\bs G_{\bs E}(\bs r, \bs r')$ in layered media. We note that the initial box in the FMM for the computation of a reaction component is set to be the smallest box containing all the equivalent polarization coordinates and effective locations. 

\subsection{Improving the efficiency of M2L evaluations}\label{sect_m2l}
Despite centers having fixed relative locations in a hierarchical box structure of FMM and thus allowing tabulation, the M2L translation typically takes major computational cost in FMM, due to its $\mathcal{O}(p^4)$ complexity, and the double integrals $\mathcal{I}_{nm\kappa}^{\nu \mu}(\bs r,\bs r',\sigma)$ from \cref{general_integral_I} that need to be numerically computed for the reaction fields.
In this section, we improve the implementation of M2L from various aspects.

\subsubsection{The number of S2T direct computations and M2L tabulations}
\begin{figure}[!ht]  
    \centering
    \includegraphics[width=0.7\linewidth]{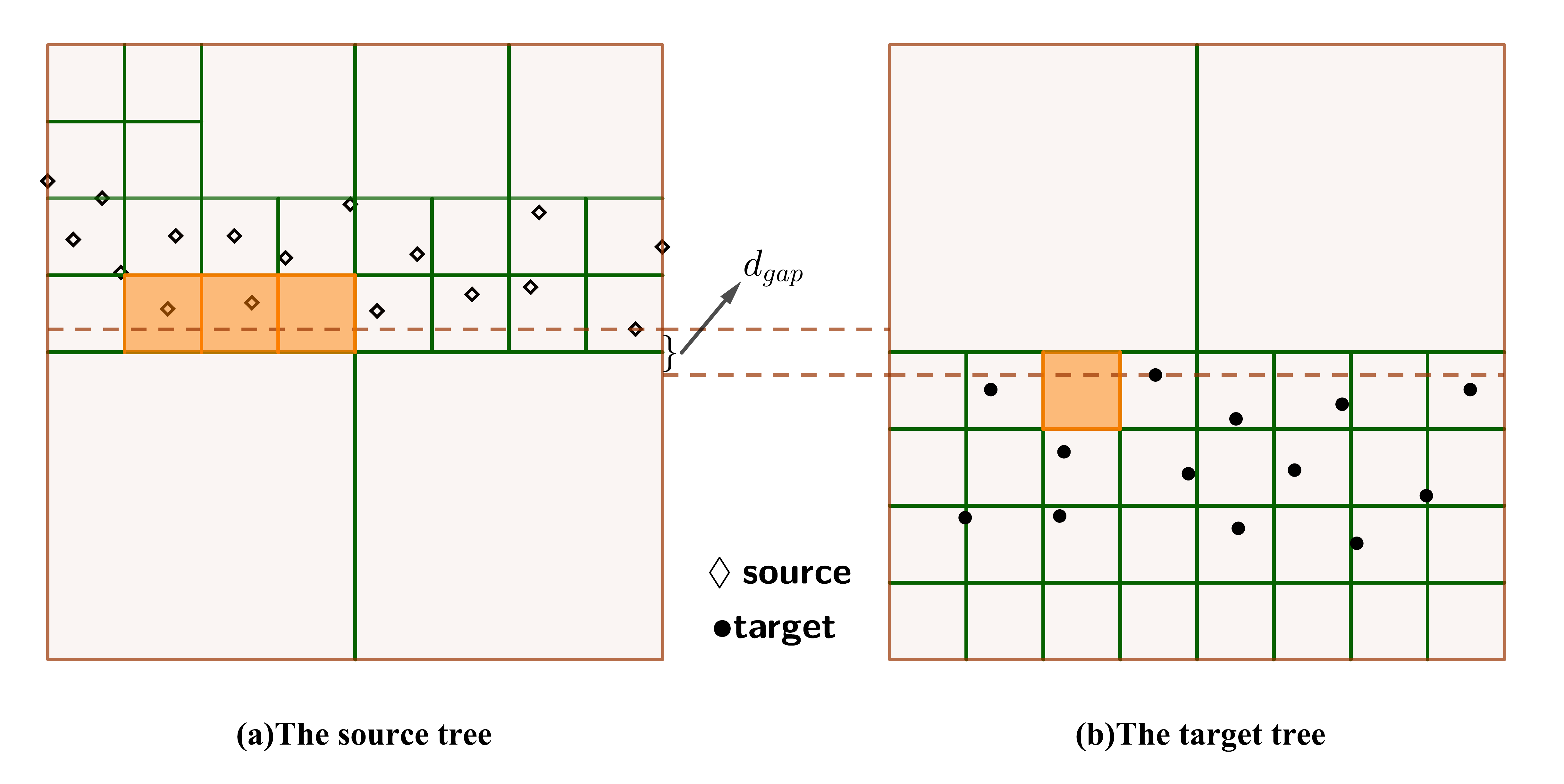}
    \caption{Cross section of an example adaptive hierarchical tree structure with adjacent non-empty source and target boxes along the interface.}
    \label{s2t_neighboring}
\end{figure}

\begin{figure}[!ht]  
    \centering
    \vspace{-15pt}
    \includegraphics[width=0.7\linewidth]{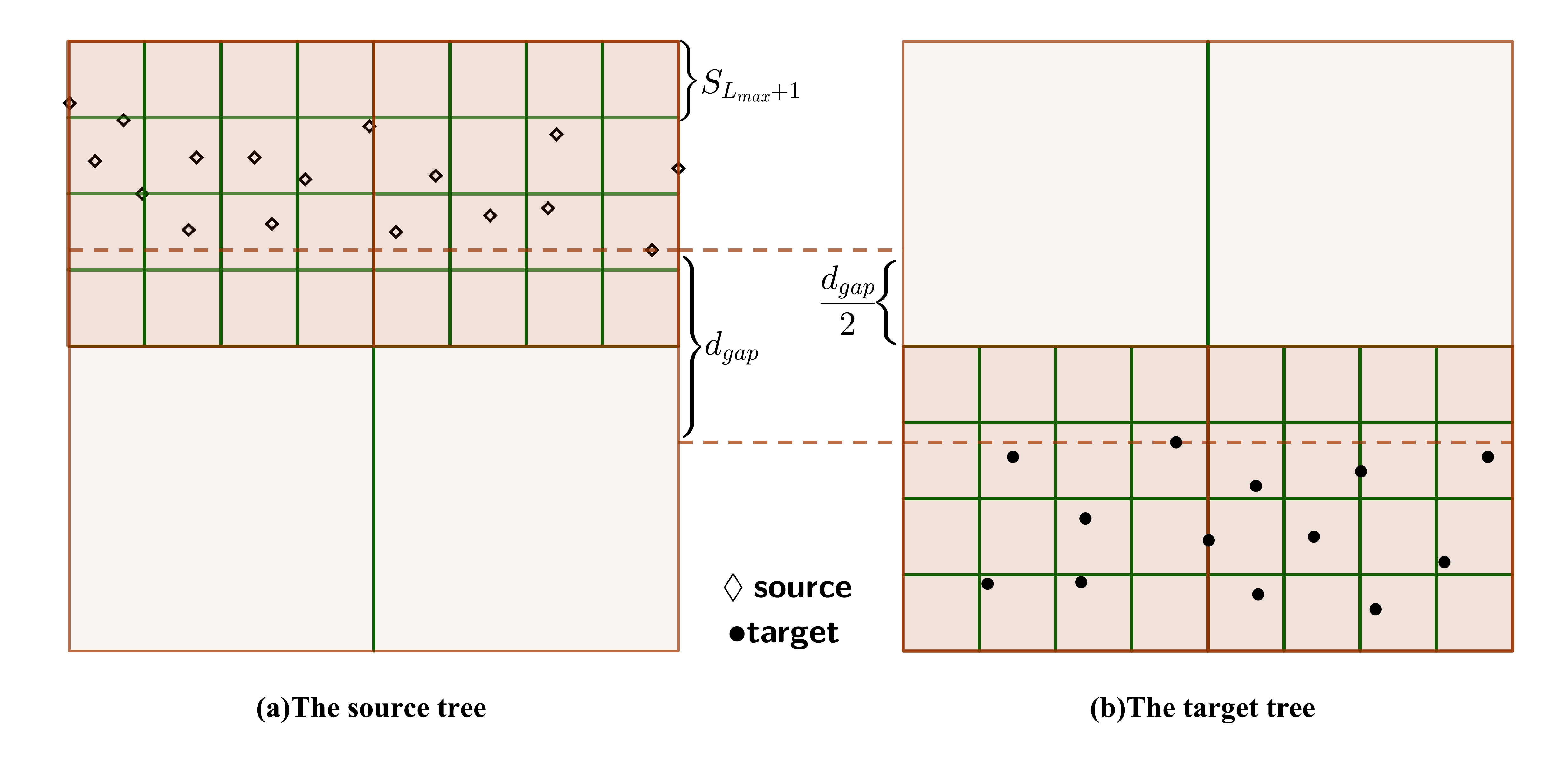}
    \caption{Cross section of an example adaptive hierarchical tree structure in which $L_{\rm max}=2$.}
    \label{m2l_empty}
\end{figure}

In the FMMs for reaction field components, only interactions between particles in adjacent boxes are computed directly. However, most of the leaf target boxes do not have neighboring source boxes as the equivalent polarization coordinates and effective locations are always separated by at least one interface. Specifically, let $d_{gap}$ be the minimum distance between equivalent polarization coordinates and the effective locations along $z$-direction. If all the size of non-empty leaf source and target boxes are smaller than $d_{gap}/2$ (see \cref{m2l_empty} for an example), then the hierarchical tree structure of FMM contains no adjacent non-empty source and target boxes.  As a result, no local direct interactions are required in this scenario. Therefore, the local direct interactions will be absent in the computation of most of the reaction field components induced by sources and targets in nonadjacent layers. For the reaction field components due to particles in the same layer or neighboring layers, we could have adjacent source and target boxes along the interface (see \cref{s2t_neighboring} for an example) where the number is generally $\mathcal O(N^{\frac{2}{3}})$, assuming $\mathcal{O}(N)$ particles are uniformly distributed in each layer. Consequently, the total number of Sommerfeld-type integrals for local direct interactions is scaled as $\mathcal O(N^{\frac{2}{3}})$. 

Within the implementation of FMM, the M2L matrices can be pre-computed for all possible $(\rho, z, z')$ determined by the target box and all source boxes in its interaction list. Again, as the equivalent polarization coordinates and effective locations are always separated by at least one interface, we only have no more than 37 different $(\rho, z, z')$ cases for all M2L in a fixed level of the tree structure. Moreover, the necessity of M2L tabulation at the current tree level $l_c$ is determined by the box size $S_{l_c}$ and $d_{gap}$, namely, when $d_{gap}\leq 2S_{l_c}$, see \cref{m2l_empty} for an example. Therefore, the maximum number of tree levels requiring M2L tabulation is given by
\begin{equation*}
	L_{\max}=\left[\log_{2}\dfrac{S_0}{d_{gap}}\right]+1,
\end{equation*}
where $S_0$ is the size of the largest box from which we start to build the tree structure. Apparently, $L_{\max}$ is small in the computation of most cases of reaction field components, as $d_{gap}$ is large when the particles are not in adjacent layers.

\subsubsection{Reducing double integrals to Sommerfeld-type integrals}
In the M2L translations and local direct interactions, there are double integrals $\mathcal{I}_{nm\kappa}^{\nu \mu}(\bs r,\bs r',\sigma)$ from \cref{general_integral_I} that need to be numerically computed.
By applying the identity 
\begin{equation*}\label{besseljidentityv}
    J_n(z)=\frac{1}{2\pi \ri^n}\int_0^{2\pi}e^{\ri z\cos\theta+\ri n\theta}d\theta,
\end{equation*}
we get
\begin{equation}\label{general_integral_I_Alternative}
    \mathcal{I}_{nm\kappa}^{\nu \mu}(\bs r,\bs r',\sigma)=(\i e^{\i\varphi})^{m-\mu+\kappa}c_{nm}^{\nu\mu}\widehat{\mathcal I}_{nm\kappa}^{\nu\mu}(\rho, z, z';\sigma), \quad \kappa = 0, 1, 2,
\end{equation}
where
\begin{equation}\label{onedimensionintegral}
\begin{split}
 \widehat{\mathcal I}_{nm\kappa}^{\nu\mu}(\rho, z, z';\sigma)=\int_0^{\infty}k_{\rho}J_{m-\mu+\kappa}(k_{\rho}\rho)\mathcal{Z}(k_{\rho},z,z') \widehat{P}_n^m\Big(\dfrac{k_z}{k}\Big)\widehat{P}_{\nu}^{\mu}\Big(\dfrac{k'_z}{k'}\Big)\sigma(k_{\rho}) dk_{\rho}
\end{split}
\end{equation}
are Sommerfeld-type integrals and $(\rho,\varphi)$ is the polar coordinates of $(x-x', y-y')$.
The cases of $\kappa=-1$ and $\kappa=-2$ can be reduced to \cref{general_integral_I_Alternative} using the identity $J_{-n}(z)=(-1)^nJ_n(z)$. Namely,
\begin{equation}\label{integralidentity}
    \mathcal{I}_{n,-m,\kappa}^{\nu,-\mu}(\bs r,\bs r',\sigma)	=(-1)^{m+\mu}e^{-2\ri (m-\mu-\kappa)\varphi}\mathcal{I}_{n,m,-\kappa}^{\nu \mu}(\bs r,\bs r',\sigma),\quad \kappa=-1, -2.
\end{equation}

\subsubsection{Reducing the number of numerical integrations in each M2L tabulation}
According to the M2L translation formula \cref{M2L_scale_with_Tnmu}, it is evident that the computation of one M2L matrix requires evaluating Sommerfeld-type integrals $\widehat{\mathcal{I}}_{nm\kappa}^{\nu \mu}$ defined in \cref{onedimensionintegral} for all $n, \nu=0, 1, \cdots, p$, $|m|\leq n$, $|\mu|\leq \nu$ and $|\kappa|\leq 2$, with a total number of $\mathcal{O}(p^4)$ Sommerfeld-type integrals. This high cost of M2L tabulation is essentially due to the $4$-entry indices $n,m,\nu,\mu$ within the product of associated Legendre functions. To simplify the M2L tabulation, we introduce a Chebyshev polynomial expansion for the products of the associated Legendre function to reduce the number of Sommerfeld-type integrals evaluated.

Let $\beta, \beta' \in \mathbb{C}$ such that
\begin{equation*}
	\sin\beta=\dfrac{k_{\rho}}{k},\quad \cos\beta = \dfrac{k_{z}}{k}, \quad \sin\beta'=\dfrac{k_{\rho}}{k'}, \quad \cos\beta' = \dfrac{k'_z}{k}.
\end{equation*}
Then,
\begin{equation*}
	\cos2\beta=1-2\frac{k_{\rho}^2}{k^2},\quad \cos2\beta'=1-2\frac{k_{\rho}^2}{k'^2}=1-\gamma^2+\gamma^2\cos2\beta,
\end{equation*}
where $\gamma=k/k'$.
Denote the product of associated Legendre functions by
\begin{equation}
	Q_{nm}^{\nu\mu}(k_{\rho})=\widehat{P}_n^m\Big(\dfrac{k_z}{k}\Big)\widehat{P}_{\nu}^{\mu}\Big(\dfrac{k'_z}{k'}\Big)=\widehat{P}_n^m\left(\cos\beta\right)\widehat{P}_{\nu}^{\mu}\left(\cos\beta'\right).
\end{equation}
We separately treat the cases $\gamma \le 1$ and $\gamma > 1$ for numerical concern in the iterative implementations. When $\gamma \le 1$, by the expansion \cref{transformedpoly} for Chebyshev polynomials,
\begin{equation}\label{transformedpolysin}
	T_j(\cos2\beta')=\sum\limits_{s=0}^jC_{js}^{ab}T_s(\cos2\beta),\quad j\geq 0\ ,\ a=\gamma^2\ ,\ b=1-\gamma^2,
\end{equation}
where $C_{js}^{ab}$ are coefficients calculated by the recurrence formula \cref{chebyshev_coefficients}.
Applying the Chebyshev polynomial expansions of the associated Legendre functions and the representation \cref{transformedpolysin}, we obtain
\begin{equation*}
\begin{split}
    Q_{nm}^{\nu\mu}(k_{\rho})&=s_{m\mu}\tau_{n,|m|}(\beta)\tau_{\nu,|\mu|}(\beta')\sum\limits_{i=0}^{K_n}\sum\limits_{j=0}^{K_{\nu}}B_{n,|m|}^iB_{\nu,|\mu|}^jT_i(\cos2\beta)T_j(\cos2\beta')\\
    &=s_{m\mu}\tau_{n,|m|}(\beta)\tau_{\nu,|\mu|}(\beta')\sum\limits_{i=0}^{K_n}\sum\limits_{j=0}^{K_{\nu}}\sum_{l=0}^jB_{n,|m|}^iB_{\nu,|\mu|}^jC_{jl}^{ab}T_i(\cos2\beta)T_l(\cos2\beta)
\end{split}
\end{equation*}
where the $B$-coefficients can be found in \cref{eq_assoc_legendre_B}, and
\begin{equation}\label{eq_sgn}
    \mathrm{sgn}(x) = (-1)^{\min(0, x)}, \quad s_{m\mu}={\rm sgn}(m){\rm sgn}(\mu).
\end{equation}
By using the identity $2T_i(x)T_j(x)=T_{i+j}(x)+T_{|i-j|}(x)$ to reduce the products of Chebyshev polynomials to sums, we arrive at an expansion of $Q_{nm}^{\nu\mu}(\kr)$ using Chebyshev polynomials
\begin{equation}\label{legendre_Q_v}
    \begin{split}
        Q_{nm}^{\nu\mu}{}&{}(\kr)
        =\frac{s_{m\mu}}{2}\tau_{n,|m|}(\beta)\tau_{\nu,|\mu|}(\beta')\sum\limits_{i=0}^{K_{n}}\sum\limits_{j=0}^{K_{\nu}}B_{n,|m|}^iB_{\nu,|\mu|}^j\times\\
        {}&{}\Big(\sum\limits_{l=i}^{j+i}C_{j,l-i}^{ab}T_{l}(\cos2\beta)+\sum\limits_{l=l_0}^{i}C_{j,i-l}^{ab}T_{l}(\cos2\beta)+\sum\limits_{l=1}^{j-i}C_{j,l+i}^{ab}T_{l}(\cos2\beta)\Big),
    \end{split}
\end{equation}
where $l_0=i-\min(i,j)$. Similarly, for $k > k'$, we have
\begin{equation}\label{legendre_Q_w}
	\begin{split}
		Q_{nm}^{\nu\mu}{}&{}(\kr)
		=\frac{s_{m\mu}}{2}\tau_{n,|m|}(\beta)\tau_{\nu,|\mu|}(\beta')\sum\limits_{i=0}^{K_{n}}\sum\limits_{j=0}^{K_{\nu}}B_{n,|m|}^iB_{\nu,|\mu|}^j\times\\
		{}&{}\Big(\sum\limits_{l=j}^{j+i}C_{i,l-j}^{ab}T_{l}(\cos2\beta')+\sum\limits_{l=l_0}^{j}C_{i,j-l}^{ab}T_{l}(\cos2\beta')+\sum\limits_{l=1}^{i-j}C_{i,l+j}^{ab}T_{l}(\cos2\beta')\Big).
	\end{split}
\end{equation}

For the Sommerfeld-type integrals $\widehat{\mathcal I}_{nm\kappa}^{\nu\mu}(\rho, z, z';\sigma)$ defined in \cref{onedimensionintegral}, let
\begin{equation}\label{chebyshev_integral_I}
	\breve{\mathcal I}_{nm\kappa}({\rho}, z, z';\tau)=\int_0^{\infty}k_{\rho}J_{m+\kappa}(k_{\rho}\rho)\mathcal{Z}(k_{\rho}, z, z')T_n(\cos 2\beta_0)\sigma(k_{\rho})\tau(\beta,\beta') dk_{\rho},
\end{equation}
where $\beta_0$ is the one of $\{\beta,\beta'\}$ satisfying $\sin \beta_0 = \kr / \min(k, k')$.
Substituting the expansions \cref{legendre_Q_v,legendre_Q_w} into \cref{onedimensionintegral}, we obtain
\begin{equation*}
	\begin{aligned}
        \widehat{\mathcal{I}}_{nm\kappa}^{\nu \mu}(\rho, z,z';\sigma)&=\frac{s_{m\mu}}{2}\sum\limits_{i=0}^{K_{n}}\sum\limits_{j=0}^{K_{\nu}}B_{n,|m|}^iB_{\nu,|\mu|}^j\mathcal A_{ij}^{m\mu\kappa}({\rho}, z, z';\tau_{n,|m|}\tau_{\nu,|\mu|})\ ,\ k\le k',\\
        \widehat{\mathcal{I}}_{nm\kappa}^{\nu \mu}(\rho, z, z';\sigma)&=\frac{s_{m\mu}}{2}\sum\limits_{i=0}^{K_{\nu}}\sum\limits_{j=0}^{K_{n}}B_{n,|m|}^jB_{\nu,|\mu|}^i\mathcal A_{ij}^{m\mu\kappa}({\rho}, z, z';\tau_{n,|m|}\tau_{\nu,|\mu|})\ ,\ k>k',
    \end{aligned}
\end{equation*}
where
\begin{equation*}
\begin{aligned}
	\mathcal A_{ij}^{m\mu\kappa}({\rho}, z, z';\tau)=&\sum\limits_{l=i}^{j+i}C_{j,l-i}^{ab}\breve{\mathcal I}_{lm\kappa}({\rho}, z, z';\tau)+\sum\limits_{l=l_0}^{i}C_{j,i-l}^{ab}\breve{\mathcal I}_{lm\kappa}({\rho}, z, z';\tau)\\
	&+\sum\limits_{l=1}^{j-i}C_{j,l+i}^{ab}\breve{\mathcal I}_{lm\kappa}({\rho}, z, z';\tau).
\end{aligned}
\end{equation*}
Therefore, the integrals \cref{onedimensionintegral} has been represented as combinations of integrals involving Chebyshev polynomials. As a result, the number of Sommerfeld-type integrals  in the reaction field M2L tabulation is reduced from $\mathcal{O}(p^4)$ to $\mathcal{O}(p^2)$.

\subsubsection{Contour deformation of Sommerfeld-type integrals}
\begin{figure}[!ht]  
    \centering
    \includegraphics[width=0.4\linewidth]{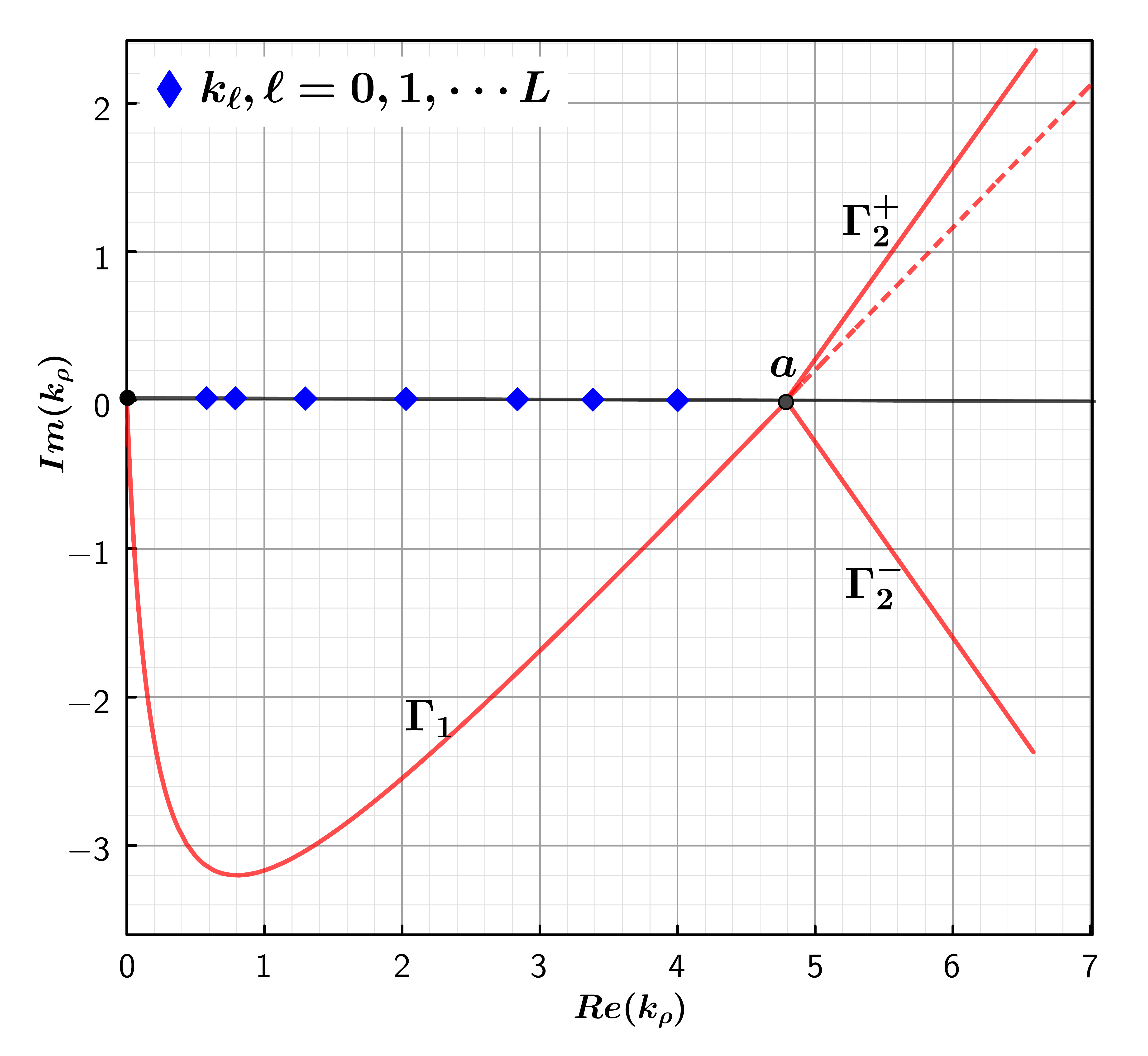}
    \caption{The deformed contours for the numerical evaluation of Sommerfeld-type integrals.}
    \label{integral_road_krho}
\end{figure}

The numerical calculation of the integral  \cref{chebyshev_integral_I} is still a challenging problem. In general, we consider integrals of the form
\begin{equation*}
    \mathcal S_{mn}(\rho, z, z')=\int_0^{\infty}J_m(k_{\rho}\rho)T_n\Big(1-2\frac{k_{\rho}^2}{k^2}\Big)f(k_{\rho}, z, z')dk_{\rho},
\end{equation*}
where $f(k_{\rho}, z, z')$ decay exponentially as $\mathfrak{Re}k_{\rho}\rightarrow\infty$. The difficulties on the computation of the integral are three folds: i) $J_n(k_{\rho}\rho)$ is highly oscillatory when $n$ or $\rho$ is large; ii) $f(k_{\rho}, z, z')$ has poles and branch cuts in the first quadrant of the complex plane of $\kr$; iii) $f(k_{\rho}, z, z')$ sometimes decays slowly, e.g. when the vertical transmission distance is small. To bypass the poles, we deform the contour for $k_{\rho}\in [0, a]$ as follows
\begin{equation*}
\begin{split}
\mathcal S_{mn}(\rho, z, z')=&\int_{\Gamma_1\cup[a,\infty)}J_m(k_{\rho}\rho)T_n\Big(1-2\frac{k_{\rho}^2}{k^2}\Big)f(k_{\rho}, z, z')dk_{\rho}
\end{split}
\end{equation*}
where $a>\max\limits_{0\leq\ell\leq L}\{k_{\ell}\}$ is a given point on the real axis such that all poles and branch cuts are on left of the line $\{k_{\rho}: \mathfrak{Re}[k_{\rho}]=a\}$, and
\begin{equation}
	\Gamma_1=\{k_{\rho}=\sqrt{2k'a\hat{u}t-a^2\hat{u}^2t^2}:0\le t<1\},\quad \hat{u}=\dfrac{k'}{a}-\ri\sqrt{1-\left(\dfrac{k'}{a}\right)^2}
\end{equation}
is the contour away from the poles and branch cuts (see \cref{integral_road_krho}). For the integral from $a$ to infinity, we use decomposition $2J_n(z)=H^{(1)}_n(z)+H^{(2)}_n(z)$ and then split the integral into two, along contours
\begin{equation}
\Gamma_2^{\pm}=\{k_{\rho}=(\Delta z\pm \ri \rho )\frac{t}{r}+a,t\ge0\},
\end{equation}
respectively. Namely,
\begin{equation*}
\begin{split}
{}&{}\int_a^{\infty}J_m(k_{\rho}\rho)T_n\Big(1-2\frac{k_{\rho}^2}{k^2}\Big)f(k_{\rho}, z, z')dk_{\rho}\\
={}&{}\frac{1}{2}\int_{\Gamma_2^+}H_m^{(1)}(k_{\rho}\rho)T_n\Big(1-2\frac{k_{\rho}^2}{k^2}\Big)f(k_{\rho}, z, z')dk_{\rho} \\
{}&{}+\frac{1}{2}\int_{\Gamma_2^-}H_m^{(2)}(k_{\rho}\rho)T_n\Big(1-2\frac{k_{\rho}^2}{k^2}\Big)f(k_{\rho}, z, z')dk_{\rho}.
\end{split}
\end{equation*}
The asymptotic behavior of the Bessel functions ensures the contour deformation and the rapid decay of the integrand as $\mathfrak{Re} k_{\rho}\rightarrow\infty$ along the new contours. Finally, a self-precision control double-exponential quadrature rule \cite{doubleexponential} is adopted to calculate the integrals along the deformed contours numerically.

\section{Numerical tests} \label{sect-num}
In this section, we present numerical results to demo-nstrate the performance of the proposed FMM for wave sources in layered media. The numerical tests are performed on a workstation, using \emph{one} CPU core of an Intel Xeon E5-2699 v4 processor (2.2 GHz) and 62 GB RAM, with GCC 12.3 compiler.

Numerical tests are conducted for both two-layers and three-layers media. Specifically, the interfaces are located at $d_0 = 0$ for the two-layers case, and at $d_0 = 0$ and $d_1 = -1.5$ for the three-layers case. The angular frequency and magnetic permeability are set as $\omega = 2.0$ and $\mu_\ell = 1.0$ for $\ell = 0, 1, 2$. The dielectric permittivities are given by $\epsilon_0 = 1.2$, $\epsilon_1 = 0.8$ for the two-layers case, and $\epsilon_0 = 1.2$, $\epsilon_1 = 0.8$, $\epsilon_2 = 1.3$ for the three-layers case. Particles are placed inside cubes of side length 1 centered at $(0.5, 0.5, 0.75)$ and $(0.5, 0.5, -0.75)$ for the two-layers configuration, and at $(0.5, 0.5, 0.75)$, $(0.5, 0.5, -0.75)$, and $(0.5, 0.5, -2.25)$ for the three-layers configuration.

{\bf Accuracy test:} We first use an example with particles uniformly distributed in the cubic domains for the accuracy test.
Let $\widetilde{\bs \Phi}_{\ell}(\bs r_{\ell i})$ be the approximated values of $\bs \Phi_{\ell}(\bs r_{\ell i})$ calculated by FMM.
We put $N = 1000$ particles randomly in each box and define $L^2$-error and maximum error as
\begin{figure}[ht!]  
		\centering
		\begin{subfigure}{0.43\textwidth}
			\includegraphics[width=\linewidth]{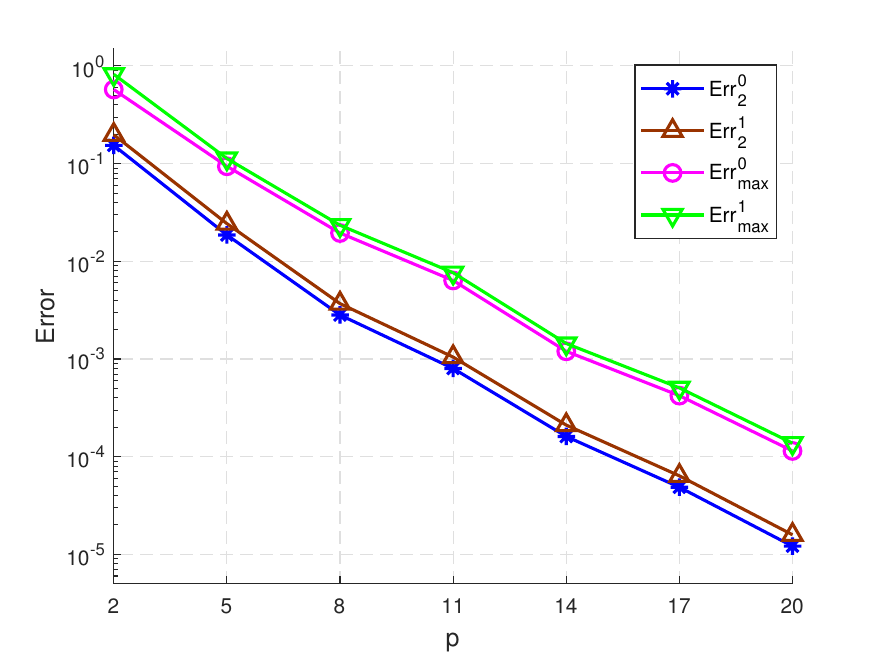}
			\caption{free-space components}
		\end{subfigure}
		\qquad
		\begin{subfigure}{0.43\textwidth}
			\includegraphics[width=\linewidth]{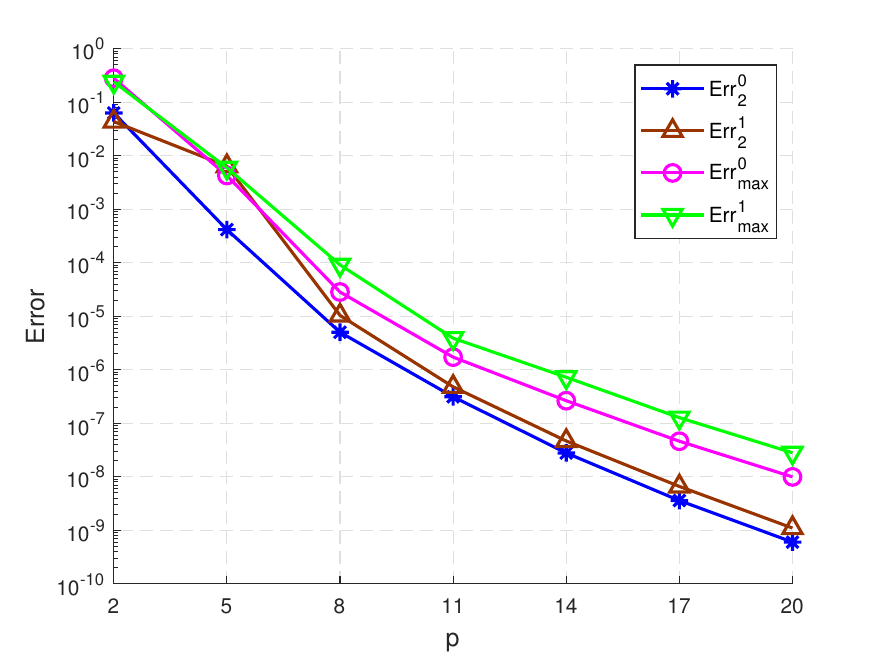}
			\caption{reaction field components}
		\end{subfigure}
		\caption{Error vs. truncation parameter $p$ in two-layer medium.}
		\label{convergence_rates_for_two-layer_vs.p}
	\end{figure}
  \begin{figure}[ht!]  
		\centering
		\begin{subfigure}{0.43\textwidth}
			\includegraphics[width=\linewidth]{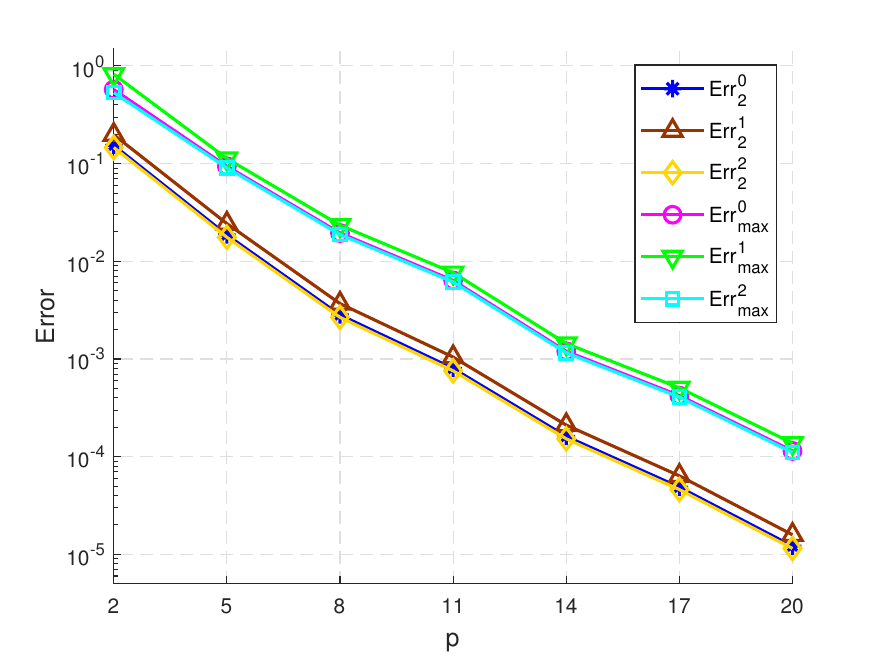}
			\caption{free-space components}
		\end{subfigure}
		\qquad
		\begin{subfigure}{0.43\textwidth}
			\includegraphics[width=\linewidth]{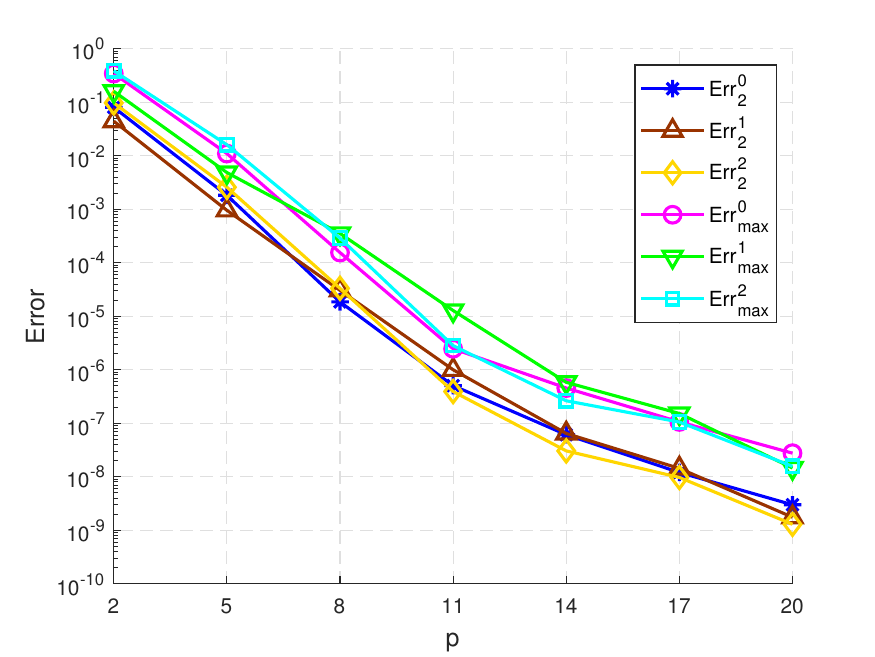}
			\caption{reaction field components}
		\end{subfigure}
		\caption{Error vs. truncation parameter $p$ in three-layer medium.}
		\label{convergence_rates_for_three-layer_vs.p}
	\end{figure}  
\begin{equation*}
    \mathrm{Err}_{2}^{\ell} = \sqrt{\dfrac{\sum_{j=1}^N \left\|\widetilde{\bs \Phi}_{\ell}(\bs r_{\ell i}) - \bs \Phi_{\ell}(\bs r_{\ell i})\right\|_2^2}{\sum_{j=1}^N \left\|\bs \Phi_{\ell}(\bs r_{\ell i})\right\|_2^2}},\quad
    \mathrm{Err}_{max}^{\ell} = \max_{1\le i\le N_{\ell}} \dfrac{\left\|\widetilde{\bs \Phi}_{\ell}(\bs r_{\ell i}) - \bs \Phi_{\ell}(\bs r_{\ell i})\right\|_2}{\left\|\bs \Phi_{\ell}(\bs r_{\ell i})\right\|_2},
\end{equation*}
for $\ell=0,1,2,\cdots,L$, where $\|\cdot\|_2$ represents the $L^2$-norm of vectors in $\mathbb R^3$.
The values of $\bs \Phi_{\ell}(\bs r_{\ell i})$ for comparison of the proposed FMM are evaluated by direct calculation of S2T interactions, using the numerical integration methods discussed in \Cref{sect_m2l}.
The convergence of the free-space components and the reaction fields as the truncation parameter $p$ increases is shown in \cref{convergence_rates_for_two-layer_vs.p} and \cref{convergence_rates_for_three-layer_vs.p}, respectively. Spectral convergence is clearly observed for both the free-space and reaction field components in the FMM. However, due to the presence of second-order derivatives in the approximations \cref{ME_Maxwell_GE,LE_Maxwell_GE}, the FMM for the free-space components typically loses about three decimal digits of accuracy compared to that for the reaction field components.

{\bf Efficiency test:} For efficiency test, we uniformly generate $M$ particles in each layer, with a total number $(L+1)M$.
The CPU time for the computation of free-space interactions ${\bs \Phi}^{f}_{\ell'}(\bs r_{\ell' i})$ and the reaction fields ${\bs\Phi}^{r}_{\ell'}(\bs r_{\ell' i})=\sum\limits_{\ell=0}^{L}\sum\limits_{\ast,\star=\uparrow,\downarrow}\bs\Phi_{\ell\ell'}^{\ast\star}(\bs r_{\ell' i}), \ell'=0,1, \cdots, L$ for $L=1, 2$ are compared in \cref{Cputime_N_p8}. 
It shows that the FMM algorithms have an $\mathcal{O}(N\log N)$ complexity while the CPU time for
the computation of reaction field components has a much smaller linear scaling constant
due to the fact that most of the equivalent polarization coordinates of source and target particles are well separated.

\begin{table}[hb!]
    \centering
    \begin{tabular}{|c|c|c|c|c|c|}
        \hline
        \# layers&\diagbox{FMM}{$M$} & 1,000 & 10,000 & 100,000 & 1,000,000 \\[1.2ex] \hline 
        \multicolumn{1}{|c|}{\multirow{4}*{2}}
        &$\Phi^{f}_{0}$ & 3.25881 & 25.1326 & 255.93 & 2337.71 \\
        \cline{2-6}
        &$\Phi^{f}_{1}$ & 3.37038 & 25.24 & 253.954 & 2432.47 \\
        \cline{2-6}
        &$\Phi_{0}^{r}$ & 3.77111 & 7.18326 & 10.58687 & 47.3383 \\
        \cline{2-6}
        &$\Phi_{1}^{r}$ & 3.77267 & 7.21007 & 10.76993 & 47.0534 \\\hline
        \multicolumn{1}{|c|}{\multirow{6}*{3}}
        &$\Phi^{f}_{0}$ & 3.4056 & 26.5547 & 274.51 & 2527.65 \\
        \cline{2-6}
        &$\Phi^{f}_{1}$ & 3.31971 & 25.6375 & 277.268 & 2517.28 \\
        \cline{2-6}
        &$\Phi^{f}_{2}$ & 3.38998 & 25.8973 & 273.253 & 2512.94 \\
        \cline{2-6}
        &$\Phi^{r}_{0}$ & 3.486836 & 7.457221 & 17.83635 & 104.475 \\
        \cline{2-6}
        &$\Phi^{r}_{1}$ & 7.281907 & 14.840553 & 35.82815 & 208.9857 \\
        \cline{2-6}
        &$\Phi^{r}_{2}$ & 3.689676 & 7.61582 & 17.96594 & 104.4376 \\\hline
    \end{tabular}
    \caption{CPU time (seconds) in two-layer and three-layer media with $p = 8$.}
    \label{Cputime_N_p8}
\end{table}

\section{Conclusion} \label{sect-conclusion}
We have proposed a fast multipole method (FMM) for the dyadic Green's functions of Maxwell's equations in layered media. By introducing equivalent polarization coordinates for sources and effective locations for targets, we addressed key challenges in far-field approximation of reaction field components using MEs and LEs. Moreover, we propose a Chebyshev polynomial expansion technique for efficient calculation of the M2L translation matrix. Numerical results confirm the $\mathcal{O}(N\log N)$ complexity and spectral accuracy of the FMM. The overall cost remains comparable to free-space FMMs, making the method practical for stratified electromagnetic problems with an arbitrary number of layers.

\appendix

\section{Transform of Chebyshev polynomials}\label{chebyshev}
For constants $a$ and $b$, the polynomial $T_j(ax+b)$ can be expanded using Chebyshev polynomials
\begin{equation}\label{transformedpoly}
	T_j(ax+b)=\sum\limits_{s=0}^jC_{js}^{ab}T_s(x),\quad j\geq 0,
\end{equation}
where the coefficients $C_{js}^{ab}$ can be evaluated using the following recurrence relations
\begin{equation}\label{chebyshev_coefficients}
    \begin{aligned}
        &	C_{00}^{ab}=1, \enspace C_{10}^{ab}=b,\enspace C_{21}^{ab}=2aC_{10}^{ab}+2bC_{11}^{ab}=4ab,\\
        &	C_{n+1,n+1}^{ab}=aC_{n,n}^{ab}, \enspace n\geq 1,\quad 
        	C_{n+1,n}^{ab}=aC_{n,n-1}^{ab}+2bC_{n,n}^{ab}, \enspace n\ge 2,\\
        &	C_{n+1,s}^{ab}=a(C_{n,s-1}^{ab}+C_{n,s+1}^{ab})+2bC_{n,s}^{ab}-C_{n-1,s}^{ab},\enspace s=2, 3,\cdots, n-1,\enspace n\ge 3,\\
        &	C_{n+1,1}^{ab}=2aC_{n,0}^{ab}+aC_{n,2}^{ab}+2bC_{n,1}-C_{n-1,1}^{ab},\enspace n\ge 2,\\
        &	C_{n+1,0}^{ab}=aC_{n,1}^{ab}+2bC_{n,0}^{ab}-C_{n-1,0}^{ab},\enspace n\ge1.
    \end{aligned}
\end{equation}

\section{Chebyshev polynomial expansions of associated Legendre functions}\label{associate_legendre}
By the trigonometric expansion of associated Legendre functions \cite{sphericalWB2018},
\begin{equation}
	\widehat P_n^m(\cos\theta)=\tau_{nm}(\theta)\sum\limits_{k=0}^{K_n}B_{nm}^kT_k(\cos2\theta), \quad m\ge0
\end{equation}
where
\begin{equation*}
    \tau_{nm}(\theta)=\begin{cases}
        1, & n\;{\rm is}\;\hbox{even},\  m\;{\rm is}\;\hbox{even},\\
        \sec\theta, & n\;{\rm is}\;\hbox{odd},\  m\;{\rm is}\;\hbox{even},\\
        \sin\theta, & n\;{\rm is}\;\hbox{odd},\  m\;{\rm is}\;\hbox{odd},\\
        \tan\theta, & n\;{\rm is}\;\hbox{even},\  m\;{\rm is}\;\hbox{odd},
    \end{cases}\quad
    K_{n}=\begin{cases}
	\dfrac{n}{2}, & n\;{\rm is}\;\hbox{even},\\
        \dfrac{n+1}{2}, & n\;{\rm is}\;\hbox{odd},\  m\;{\rm is}\;\hbox{even},\\
 	\dfrac{n-1}{2}, & n\;{\rm is}\;\hbox{odd},\  m\;{\rm is}\;\hbox{odd},
    \end{cases}
\end{equation*}
and
\begin{equation}\label{eq_assoc_legendre_B}
B_{nm}^{k} = \begin{cases}
A_{nm}^k, & n \text{ is even, }m \text{ is even}, \\
\dfrac{1}{2m}\left(d_{n,m+1}^{+}B_{n,m+1}^k+d_{n,m-1}^{-}B_{n,m-1}^k\right), & n \text{ is even, }m \text{ is odd}, \\
\dfrac{1}{2m}\left(f_{nm}^{-}B_{n-1,m+1}^k+f_{nm}^{+}B_{n-1,m-1}^k\right), & n \text{ is odd, }m \text{ is odd, }m < n, \\
\dfrac{1}{2m}f_{nm}^{+}B_{n-1,m-1}^k & n \text{ is odd, }m \text{ is odd, }m = n, \\
a_{nm}B_{n-1,m}^k+a_{n+1,m}B_{n+1,m}^k & n \text{ is odd, }m \text{ is even},
\end{cases}
\end{equation}
and $A_{nm}^{k}$ has backward recursion on $k$
\begin{equation*}
    A_{nm}^{n/2 + 1} = 0, \quad A_{nm}^{n/2} = (-1)^{m/2} \frac{\Gamma(n+1/2)}{\pi}\sqrt{\frac{2n+1}{(n-m)!(n+m)!}},
\end{equation*}
\begin{equation*}
    A_{nm}^{k} = a_{nm}^{k+2}A_{nm}^{k+1} + b_{nm}^{k+2}A_{nm}^{k+2}, \quad k = \frac{n}{2} - 1, \frac{n}{2} - 2, \cdots, 1,
\end{equation*}
\begin{equation*}
    A_{nm}^{0} = \begin{cases}
        (a_{nm}^{2}A_{nm}^{1} + b_{nm}^{2} A_{nm}^{2})/2, & l \ge 4, \\
        \frac{n(n+1) - 2}{2n(n+1) - 4m^2} A_{nm}^{1}, & l = 2,
    \end{cases}
\end{equation*}
where the constants
\begin{equation*}
    a_{nm}^{k} = \frac{2(2m^2-n(n+1)+4(k-1)^2)}{2(k-2)(2k-3)-n(n+1)}, \quad b_{nm}^k =\frac{n(n+1)-2k(2k-1)}{2(k-2)(2k-3)-n(n+1)}.
\end{equation*}
The case $m < 0$ is dealt with the formula $\widehat P_n^{-m}(z)=(-1)^m\widehat P_n^m(z)$. In combination,
\begin{equation}\label{ALFChebyexp}
	\widehat P_n^m(\cos\theta)=\mathrm{sgn}(m)\tau_{n,|m|}(\theta)\sum\limits_{k=0}^{K_n}B_{n,|m|}^kT_k(\cos2\theta), \quad -n \le m \le n.
\end{equation}

\bibliographystyle{siamplain}
\bibliography{reference}

\end{document}